\documentclass[12pt,a4paper]{article}

\title{Conjugacy classes of diffeomorphisms of the interval in $\mc{C}^{1}$-regularity }

\usepackage[utf8]{inputenc}
\usepackage[french]{babel}
\usepackage[T1]{fontenc}

\usepackage{amsmath}
\usepackage{amsfonts}
\usepackage{amsthm}
\usepackage{amssymb}
\usepackage{makeidx}
\usepackage{mathrsfs}
\usepackage{euscript}
\usepackage{geometry}
\usepackage{graphicx}
\usepackage{textcomp}
\usepackage{stmaryrd}

\newcommand{\limite}[2]{\mathop{\longrightarrow}
\limits_{\mathrm{#1}}^{\mathrm{#2}}}
\DeclareMathOperator{\cv}{\limite{n \to \infty}{}}

\newtheorem{lem}{Lemma}
\newtheorem{coro}{Corollary}
\newtheorem{prop}{Proposition}
\newtheorem{rema}{Remark}
\newtheorem{theo}{Theorem}
\newtheorem{defi}{Definition}[section]
\newtheorem{exem}{Example}

\newtheorem{affi}{Claim}
\newtheorem{ques}{Question}
\newtheorem{conj}{Conjecture}

\newtheorem{lede}{Lemma-Definition}

\newcommand{\N}{\mathbb{N}}
\newcommand{\Z}{\mathbb{Z}}

\newcommand{\R}{\mathbb{R}}

\newcommand{\mc}[1]{\mathcal{#1}}

\newcommand{\D}{D_{\alpha,\beta}}

\newcommand{\crochetl}{\left[\kern-0.15em\left[}
\newcommand{\crochetr}{\right]\kern-0.15em\right]}

\author{\'Eglantine Farinelli}

\begin{document}

\def\restriction#1#2{\mathchoice
              {\setbox1\hbox{${\displaystyle #1}_{\scriptstyle #2}$}
              \restrictionaux{#1}{#2}}
              {\setbox1\hbox{${\textstyle #1}_{\scriptstyle #2}$}
              \restrictionaux{#1}{#2}}
              {\setbox1\hbox{${\scriptstyle #1}_{\scriptscriptstyle #2}$}
              \restrictionaux{#1}{#2}}
              {\setbox1\hbox{${\scriptscriptstyle #1}_{\scriptscriptstyle #2}$}
              \restrictionaux{#1}{#2}}}
\def\restrictionaux#1#2{{#1\,\smash{\vrule height .8\ht1 depth .85\dp1}}_{\,#2}}

\maketitle

{\small
\centerline{\bf Abstract} 

~~In this paper we consider the conjugacy classes of diffeomorphisms of the interval, endowed with the $C^1$-topology. We present several results in the spirit of the one below : 

Given two diffeomorphisms $f,g$ of the interval $[0;1]$ without hyperbolic fixed point, we give a complete answer to the two following questions : 
\begin{itemize}
\item under what conditions does there exist a sequence of  smooth conjugates $h_n f h_n^{-1}$ of $f$ tending to $g$ in the $C^1$-topology ? 

\item under what conditions does there exist a continuous path of $C^1$-diffeomorphisms $h_t$ such that $h_t f h_t^{-1}$ tends to $g$ in the $C^1$ topology ?
\end{itemize} 

 ~We present also some consequences of these results as regards the study of the $C^1$-centralizers of $C^1$-contractions of $[0;+\infty[$ ; for instance we exhibit a $C^1$-contraction whose centralizer is uncountable and abelian, but is not a flow. 

}
\renewcommand{\thefootnote}{}


\section{Introduction}

\subsection{Conjugacy classes of diffeomorphisms of the interval}

~~We consider diffeomorphisms of $[0;1]$ whose fixed points are precisely $0$ and $1$. 
In 1970, N. Kopell \cite{Ko} showed that, for an open and $\mc{C}^2$-dense set of such diffeomorphisms $f$, the only $C^1$-diffeomorphisms commuting with $f$ are the iterates of $f$, i.e. the elements of the group $\{f^i, i\in\Z\}$ : one says that $f$ has trivial centralizer. One of the keys of this theorem is that the $\mc{C}^1$-centralizers of $f$ on $[0;1[$ and on $]0;1]$ are both one parameter-groups (more precisely, they are flows) ; let us denote them respectively by $f^-_t$ and $f^+_t$. The centralizer of $f$ on the segment $[0;1]$ corresponds to the times $t$ such that $f^-_t=f^+_t$.  Comparison between these two flows is described by Mather invariant of $f$ : this invariant depends continuously of $f$ in $\mc{C}^2$-topology, and cancels when $f^-_t=f^+_t$ for all $t$. 

Mather invariant is a $\mc{C}^1$-conjugacy invariant among $\mc{C}^2$-diffeomorphisms. However, it does not vary continuously with respect to $\mc{C}^1$-perturbations of the diffeomorphism $f$ : \cite{BCVW} shows that a $\mc{C}^1$-small perturbation of $f$ cancels its Mather invariant, and thus the diffeomorphism is embedded in a flow. The argument mentioned in \cite{BCVW} seemed to indicate that an arbitrarily small $\mc{C}^1$-perturbation of $f$ could transform Mather invariant as wanted. On the other hand, Mather invariant represents, together with the values of the derivatives at the extremities -- when these are not equal to $1$ -- a complete $\mc{C}^1$-conjugacy invariant among $\mc{C}^2$-diffeomorphisms of $[0;1]$ without fixed point in the interior.

This reasoning leads to the following conjecture, which appears already in \cite{BCVW} :
\begin{conj}\label{conj}
For all $\alpha > 1$ and $0 < \beta < 1$, one considers the set $D_{\alpha,\beta}$ of
diffeomorphisms $f\colon [0;1]\to[0;1]$ whose fixed points are exactly $0$ and $1$, and such that
$Df(0) =\alpha$ and $Df(1) =\beta$.

 $\mc{C}^1$-conjugacy classes are all $\mc{C}^1$-dense in $D_{\alpha,\beta}$. In other words,
given two diffeomorphisms $f$, $g\in D_{\alpha,\beta}$, there exists a diffeomorphism, arbitrarily $\mc{C}^1$-close to $g$ which is conjugate to $f$ by a diffeomorphism of $[0;1]$.
\end{conj}

Our first result answers entirely this conjecture : conjugacy classes are dense in $D_{\alpha,\beta}$.  Actually, we show a a little stronger result, which implies that each diffeomorphism $g\in D_{\alpha,\beta}$ is \emph{accessible} from $f$ by a path in its differentiable conjugacy class. We will see that this stronger notion is important for the question of centralizers. 

More precisely : 

\begin{defi}\label{d.classeconj}
Given $f,g\in\mc{D}iff^{1}([0;1])$, one says that $f$ is \emph{isotopic to $g$ by $\mc{C}^1$-conjugacy} (or simply \emph{isotopic by conjugacy to $g$})
if there exists a $\mc{C}^{1}$-continuous path $(h_t)_{t\in[0,1[}$,  $h_{t}\in\mc{D}iff^{1}([0;1])$, such that : 
$h_{0}=id$, and
$h_{t} f h_{t}^{-1} \limite{t \to 1}{} g$. Under these conditions, 
the path $(h_{t} f h_{t}^{-1})_{t\in [0;1[}$ will be called \emph{isotopy by conjugacy from $f$ to $g$}.
\end{defi}

\begin{rema} This notion is \emph{a priori} a little more restrictive than the fact that $g$ is \emph{accessible from $f$ by a path in its conjugacy class} (i.e: $g$ is the limit of a continuous path $(h_tf h_t^{-1})_{t\in[0,1[}$ such that $h_t\in \mc{D}iff^{1}([0;1])$ for all $t$ and $h_0=id$, without requiring any continuity of $h_t$ with respect to $t$). 

\end{rema}

\begin{ques} Given $f,g$ two diffeomorphisms of $[0;1]$, under what condition might there exist an isotopy by $\mc{C}^1$-conjugacy from $f$ to $g$ ?  
\end{ques} 

We answer first this question for diffomorphisms without fixed point in $]0;1[$, solving in particular by this way Conjecture~\ref{conj} : 

\begin{theo}\label{t.4} Let $f$ and $g$ be two diffeomorphisms of $[0;1]$ without fixed point, except $0$ and $1$, and greater than the identity on $]0;1[$. Then there exists an isotopy by conjugacy from $f$ to $g$ if and only if $f$ and $g$ have same derivatives at $0$ and $1$ : $Df(0)=Dg(0)$ and $Df(1)=Dg(1)$.
\end{theo}

Continuity of  Mather invariant in $\mc{C}^2$-topology implies that Theorem \ref{t.4} is not true in $\mc{C}^r$-topology, for all $r\geq 2$.  

Actually, Theorem~\ref{t.4} is a consequence of more general Theorem~\ref{t1.transfochemin} : each continuous path of diffeomorphisms linking $f$ to $g$ without neither changing derivatives on the boundary nor introducing any fixed point in the interior can be $C^1$-approached by an isotopy by conjugacy from $f$ to $g$.  

\vskip 5mm

~~Let us consider now diffeomorphisms of $[0;1]$ without hyperbolic fixed point in $]0;1[$, but without restriction about the set of non hyperbolic fixed point. We will see that, in this context, the \emph{existence of an isotopy by conjugacy from $f$ to $g$}  is a more restrictive notion than the fact that \emph{$g$ is an accumulation point of the  conjugacy class of $f$}.



Theorem~\ref{t1.généralisation} above provides a necessary and sufficient condition for each of those properties. In its statement appears the notion of \emph{signature} of a diffeomorphism $f$ of $[0;1]$ (see definition~\ref{d.signature}) : the signature of $f$ is a pair $((C=\{C_i\}_{i\in I},\prec),\sigma)$, where :
\begin{itemize} 
\item $(C,\prec)$
 is a countable ordered set. Roughly speaking, $C$ is the set of maximal intervals on which the sign of $f$ does not change, ordered by their position in $[0;1]$ ;
\item $\sigma$ is a map from $C$ into $\{+,-\}$  which associates to each interval the sign of $f-id$ on it. 
\end{itemize} 
 


The most general résult of this paper is the following :

\begin{theo}\label{t1.généralisation}
Let $f$ and $g$ be two non-decreasing $\mc{C}^1$-diffeomorphisms of $[0;1]$ without hyperbolic fixed point except possibly $0$ and $1$, and such that $Df(0)=Dg(0)$, and $Df(1)=Dg(1)$.
Let $((C,\prec),\sigma)$ and $((C',\prec'),\sigma')$ denote the respective signatures of $f$ and $g$.
\\
Then :
\begin{enumerate}
\item There exists a conjugacy isotopy from $f$ to $g$ if and only if there exists an injective and order preserving map $\Phi : C'\longrightarrow C$ such that, for all $c'\in C'$, $\sigma'(c')=\sigma(\Phi(c'))$.
\item There exists a sequence of conjugates of $f$ converging to $g$ if and only if, for all finite subset $\gamma'$ of $C'$, there exists an injective and order preserving map $\Phi : \gamma'\longrightarrow C$ such that, for all $c'\in \gamma'$, $\sigma'(c')=\sigma(\Phi(c'))$.
\end{enumerate}
\end{theo}


\begin{exem}\label{e.classeconj} 
Let us consider $f,g\in \mc{D}iff^{1}([0;1])$ without hyperbolic fixed point and such that
$\mathrm{Fix}(f)$ is a sequence tending to $0$, and the sign of $f-id$ changes at each fixed point ; $\mathrm{Fix}(g)$ consists in a sequence tending to $0$ and a sequence tending to $1$, and the sign of $g-id$ changes at each fixed point. 
Namely, they can been chosen as follows :

\begin{itemize}
\item sign of $f-id$ is : 
\begin{itemize}
\item[$\centerdot$] $<0$ on $]\frac 1{2n};\frac 1{2n-1}[$ and 
\item[$\centerdot$] $>0$ on $]\frac 1{2n+1};\frac 1{2n}[$ for all $n>0$ ;
\end{itemize}
\item sign of $g-id$ is 
\begin{itemize} 
\item[$\centerdot$] $<0$ on $]\frac 1{2n};\frac 1{2n-1}[$ and on $]1-\frac 1{2n-1};1-\frac 1{2n}[$

and
\item[$\centerdot$] $>0$ on $]\frac 1{2n+1};\frac 1{2n}[$, on $]1-\frac 1{2n};1-\frac 1{2n+1}[$ and on $]\frac 13;\frac 23[$ for all $n>1$.
\end{itemize}
\end{itemize}
Then there exists a sequence $(h_n)_{n\in\N}$ of diffeomorphisms of $[0;1]$ such that $h_nf h_n^{-1}$ converges to $g$, but there does not exist any isotopy by conjugacy from $f$ to $g$. Indeed, one can easily check that the signature of $f$ is $-\N$ with $\sigma_f(-n)=(-1)^{n+1}$ whereas the signature of $g$ is $\Z$ with $\sigma_g(n)=(-1)^n$. 
\end{exem}

\subsection{$\mc{C}^r$-centralizer of contractions of $[0;+\infty[$}

~~The \emph{centralizer $\mc{C}^{r}(f)$}, for $r\geqslant 0$, of a $\mc{C}^{r}$-diffeomorphism $f$ of a manifold $M$ is the set
$\mc{C}^{r}(f)=\{g\in Diff^{r}(M)\; ; \;gf=fg\}$.
This set is a group and contains the group $<f>=\{f^i\}_{i\in\Z}$ generated by $f$. One says that $f$ has \emph{trivial centralizer} when $\mc{C}^{r}(f)$ is as small as possible, that is $\mc{C}^{r}(f)=<f>$.
\\

We are interested here in diffeomorphisms of the half-line $[0;+\infty[$, and more precisely in \emph{contractions}, that is diffeomorphisms having $0$ as only fixed point and which are smaller than the identity on $]0;+\infty[$. In this context, size and structure of centralizers of diffeomorphisms of $[0;+\infty[$ depend essentially on the regularity of considered diffeomorphisms.
\\
For example, each $\mc{C}^{0}$-contraction $f$ of $[0;+\infty[$ is conjugate (in restriction to $]0;+\infty[$) to the translation $x\longmapsto x-1$ of $\R$,
so that the centralizer of $f$ is conjugate to the group of lifts on $\R$ of diffeomorphisms of the circle. Consequently, the $\mc{C}^0$-centralizer of a contraction has large size and does not depend, up to a conjugacy, on the considered diffeomorphism.
\\

~~By simply working beyond the $\mc{C}^{0}$ context, very opposite situations can occur : if $f$ is an homothety of $[0;+\infty[$, then the only homeomorphisms commuting with $f$ and which are differentiable at zero are homotheties ; in particular $\mc{C}^1(f)$ is the group of homotheties. 
In the case of a $\mc{C}^{2}$-contraction $f$ of $[0;+\infty[$, the situation is quite similar : the $\mc{C}^1$-centralizer is then isomorphic to $\R$ (it is a flow). This follows from Kopell's lemma (see \cite{Ko}) and from Szekeres' theorem (see \cite{Sz}) or from Sternberg's theorem (see \cite{St}). Consequently :

\emph{The $\mc{C}^{r}$-centralizer of a contraction of $[0;+\infty[$, for $r\geqslant 2$, is isomorphic to a subgroup of $\R$}.

\vskip 2mm

That is what was stated in Kopell's lemma : if $f$ is a $\mc{C}^2$-contraction of $[0;+\infty[$ and $g$ is a $C^1$-diffeomorphism commuting with $f$, then, if $g$ has another fixed point than $0$, $g$ must be the identity. From that, one can deduce immediately that the  $\mc{C}^1$-centralizer of $f$ is an ordered and archimedean group, and thus from Hölder's theorem is a subgroup of $\R$. 
 
Result above lead us to the natural following question :
\begin{ques} Which subgroups of $\R$ appear as $\mc{C}^r$-centralizers of contractions of $[0;+\infty[$? 
\end{ques}
Let us recall about that question that F. Sergeraert, in \cite{Se}, gave explicit examples of $\mc{C}^\infty$-contractions of $[0;+\infty[$ whose $\mc{C}^{r}$-centralizers, for $r\geqslant 2$, are not whole $\R$, and also that H. Eynard, from this construction, presented in \cite{Ey} examples where the centralizer is a proper subgroup of $\R$ containing a Cantor set.

~~Between these two extreme situations, we find the case of the $\mc{C}^1$-centralizer of $\mc{C}^{1}$-contractions of $[0;+\infty[$, in which we are interested in this paper. This context is well distinct from the $\mc{C}^{2}$-context, in particular because, since Kopell's lemma is nonvalid in the $\mc{C}^{1}$-context, the centralizer presents a greater flexibility than in $\mc{C}^{2}$-context, and also because in compensation the nonvalidity of Szekeres' theorem does not ensure anymore the existence of a centralizer as big as $\R$.
\subsection{Embedding of a group in a centralizer and isotopy by conjugacy to identity}

~~In this context of $\mc{C}^1$-centralizers of $\mc{C}^{1}$-contractions, very various situations can occur ; indeed :
\begin{itemize}
\item Togawa, in \cite{T}, showed that, for all map in a $G_{\delta}$-dense of the set of $\mc{C}^1$-contractions, its centralizer is trivial. 
\item We will see as a consequence of our results that there exist contractions of $[0;+\infty[$ whose centralizer contains nontrivial diffeomorphisms with other fixed points than $0$ (going by this way against the conclusions of Kopell's lemma). It was already known  that such examples existed (see in particular \cite{FF}). 
\end{itemize}

Existence of $C^1$-counter-examples to Kopell's lemma lead to the natural question : given an interval of $]0;+\infty[$, which group of diffeomorphisms with support in this interval can be embedded in the centralizer of a contraction ? First we have to formulate this question in more precise words.

\begin{defi}
Given a segment $J\subset]0;+\infty[$, and
a group $G$ of diffeomorphisms of $J$, this group $G$ will be said 
\emph{embeddable in the centralizer of a contraction} if there exists a 
$\mc{C}^{1}$-contraction $f$ of $[0;+\infty[$, with $0$ as only 
fixed point, and a subgroup $G_{0}$ of $\mc{C}^{1}(f)$ such that
$\{\restriction{g}{J} ; g\in G_{0}\}=G $. 
\end{defi}

Our first result answers this question by characterizing groups which are embeddable in the  $\mc{C}^1$-centralizer of a contraction.
Before stating this result, we extend the notion of isotopy by conjugacy from a diffeomorphism to another (definition \ref{d.classeconj}) to the notion of isotopy by conjugacy from a group to a diffomorphism :

\begin{defi}\label{d.plongement}
A subgroup $G\subset \mc{D}iff^1([0;1])$ is said to be \emph{isotopic by conjugacy to a $\mc{C}^1$-diffeomorphism $g$} de $[0;1]$
if there exists a continuous path $(h_{t})_{t\in[0;1[}$ of $\mc{C}^{1}$-diffeomorphisms of $[0;1]$ such that, for all $f\in G$, $(h_{t}fh_{t}^{-1})_{t\in[0;1[}$ is an isotopy by conjugacy from $f$ to $g$. 
\end{defi}

\begin{theo}\label{t.1}
A group of diffeomorphisms of a segment of $]0;+\infty[$ is embeddable in the centralizer of a contraction if and only if there exists an isotopy by $\mc{C}^1$-conjugacy from this group to $id $. 
\end{theo}

Thus we want now to answer the following question : 
\begin{ques} Which groups of diffeomorphisms of $[0;1]$ are isotopic by conjugacy to the identity ? 
\end{ques}
An obvious obstruction is the existence of a hyperbolic fixed point (i.e. with derivative different than $1$) for at least one element of the group. The question is whether there exist other obstructions. 

In this paper, we answer this question in the case where $G$ is generated by an only diffeomorphism. We prove :

\begin{theo}\label{t.2} A $\mc{C}^1$-diffeomorphism $f$ of $[0;1]$ is isotopic to the identity by $\mc{C}^1$-conjugacy if and only if all its fixed points are tangent to the identity. 
\end{theo}

We think that Theorem~\ref{t.2} is not true if the isotopy is asked to be a $\mc{C}^2$-continuous path of $\mc{C}^2$-conjugates. 

As a corollary, it follows immediately that :

\begin{coro} If $J$ is a nontrivial segment of $]0;+\infty[$ and if $g_J$ is a diffeomorphism of $J$ whithout hyperbolic fixed point, then $g_J$ is embeddable in the $\mc{C}^1$-centralizer of a contraction : there exists a contraction $f$ of $[0;+\infty[$ and a diffeomorphism $g\in\mc{C}^1(f)$ whose restriction to $J$ is $g_J$. 
\end{coro}

The following result exhibits $C^1$-contractions with big centralizers :

\begin{coro}\label{t.3} Let $\mc{I}=\{I_i\}_{i\in\N}$ be a sequence of intervals of $[0;1]$ whose interiors are pairwise disjoint and, for all $i\in\N$, let $g_i$ be a diffeomorphism of $[0;1]$ with support in $I_i$, whose derivatives at fixed points are $1$. We assume that, for all $i\in\N$, $g_i$ tends to the identity in $\mc{C}^1$-topology : $\lim\limits_{i\to+\infty} \Vert g_i-id\Vert_1=0$. Then :
\begin{itemize}
\item for all bounded sequence of integers $\bar n=\{n_i\}_{i\in\N}$, the map $g_{\bar n}$ which coincides with $g_i^{n_i}$ on $I_i$ for all $i\in\N$ and with the identity outside the union of the $I_i$'s is a diffeomorphism of $[0;1]$ ;
\item the set of all $g_{\bar n}$ where $\bar n$ is a bounded sequence of integers is an abelian and uncountable group, denoted by $G_{\mc{I}}$ ;
\item the group $G_{\mc{I}}$ is isotopic by conjugacy to the identity.  
\end{itemize}
\end{coro}

\underline{Proof :} One checks that the set $G_{I}$ is a group of diffeomorphisms. From Theorem \ref{t.2}, if $\bar{n}$ is the constant sequence equal to $1$, then $g_{\bar{n}}$ is isotopic to the identity by conjugacy, by the isotopy $(h_{t}g\bar{n}h_{t}^{-1})_{t\in[0;1[}$. It follows that, for all $i$,
$(h_{t}g_{i}h_{t}^{-1})_{t\in[0;1[}$ is an isotopy by conjugacy from $g_{i}$ to the identity.
Thus, for all bounded sequence $\bar{n}$ of integers, $(h_{t}g\bar{n}h_{t}^{-1})_{t\in[0;1[}$ is an isotopy by conjugacy from $g_{\bar{n}}$ to the identity.

\hfill$\Box$\medskip

Further to a first informal version of this paper, Andr\`es Navas made us know a clever and very simple argument for proving Theorem~\ref{t.2}, which is based on the cohomological equation. His result can also be applied in the context of diffeomorphisms of the circle with irrational number of rotation. Here is what he shows : 

\begin{theo}[Navas]\label{t.navas} Each diffeomorphism of $[0;1]$ without hyperbolic fixed point is isotopic by conjugacy to the identity. Each diffeomorphism of the circle  with irrational number of rotation $\alpha$ is isotopic by conjugacy to the rotation $R_\alpha$. 
\end{theo}

\paragraph{Acknowledgements :} I would like to thank C. Bonatti for his significant contribution to the conception and drafting of this paper, as well as A. Navas for his
contribution. I would like to thank also H. Eynard, A. Wilkinson and N. Guelman for their help and interest in my work.


\section{Groups which are embeddable in a centralizer : proof of Theorem~\ref{t.1}}

~~In order to prove that theorem, we will use the following lemma :

\begin{lem}\label{l1.centralisateur et isotopie}
Let $G$ be a group of diffomorphisms of a segment $J$. There exists an isotopy from $G$ to $id$ if and only if there exists a sequence $(\varphi_{n})_{n\in\N}$ of diffeomorphisms of this segment, converging to the identity in $\mc{C}^{1}$-topology, such that $\varphi_{0}=id$ and such that, for all $g\in G$, the sequence
$\varphi_{n}\ldots\varphi_{1}\varphi_{0}g\varphi_{0}^{-1}\varphi_{1}^{-1}\ldots
\varphi_{n}^{-1}$
converges to the identity with respect to the  $\mc{C}^{1}$-topology.
\end{lem}

\underline{Proof of Lemma \ref{l1.centralisateur et isotopie} :}
$\bullet$ Let us assume that there exists an isotopy $(h_{t}gh_{t}^{-1})_{t\in[0;1[}$ 
from each element $g$ of $G$ to the identity.
\\
\underline{Claim :} There exists a sequence of strictly positive numbers $(\varepsilon_{n})_{n\in\N}$ converging to $0$ ; there exists an increasing sequence $(t_{n})_{n\in\N}$ of numbers in $[0;1[$ converging to $1$, such that $t_{0}=0$ and such that, for all $n\in\N$, one has $\Vert h_{t_{n+1}}h_{t_{n}}^{-1}-id\Vert_{1}<\varepsilon_{n}$.
\\

Proof of the claim :
\\
Given an increasing sequence $(t'_{n})_{n\geqslant 1}$ of numbers in $[0;1[$ converging to $1$, by compacity of the set $\{h_{t} ; t\in[t'_n;t'_{n+1}]\}$, for all $n\geqslant 1$ there exist a whole number $K_n$ and a sequence $(t^n_k)_{k\in \llbracket 1;K_n\rrbracket}$ of real numbers of $[0;1[$ such that 
$\Vert h_{t^n_{k+1}}h_{t^n_{k}}^{-1}-id\Vert_1<\frac{1}{n}$ for all $k\in\llbracket 1;K_n-1\rrbracket$.
By concatenating those sequences, $(t^n_k)_{k\in \llbracket 1;K_n\rrbracket}$ on one hand, and on the other hand the constant sequences equal to $\frac{1}{n}$ for $n\geqslant 1$, we obtain respectively the wanted sequences $(t_n)_{n\in\N}$ and $(\varepsilon_n)_{n\in\N}$.

\hfill$\Box$\medskip

~~Let $(\varepsilon_{n})_{n\in\N}$ and $(t_{n})_{n\in\N}$ be sequences as described in the claim above, and let us define the following $\mc{C}^{1}$-diffeomorphisms of $[0;1]$ : $\varphi_{0}=h_{t_{0}}$ and, for all $n\geqslant 1$, 
$\varphi_{n}=h_{t_{n}}h_{t_{n-1}}^{-1}$. Then the sequence $(\varphi_{n})_{n\in\N}$ converges to the identity in $\mc{C}^{1}$-topology, as well as the sequence
$(\varphi_{n}\ldots\varphi_{1}\varphi_{0}g\varphi_{0}^{-1}\varphi_{1}^{-1}\ldots
\varphi_{n}^{-1})_{n\in\N}=(h_{t_{n}}gh_{t_{n}}^{-1})_{n\in\N}$, since $h_{t}gh_{t}^{-1}$ is assumed to converge to $id$ when $t$ tends to $1$.
\\

$\bullet$ Conversely, let us assume that there exists a sequence $(\varphi_{n})_{n\in\N}$ of diffeomorphisms of $J$
converging to $id$ when $n$ tends to $+\infty$, 
such that $\varphi_{0}=id$ and such that, for all $g\in G$, $\varphi_{n}\ldots\varphi_{1}\varphi_{0}g\varphi_{0}^{-1}\varphi_{1}^{-1}\ldots
\varphi_{n}^{-1}$
converges to the identity with respect to the $\mc{C}^{1}$-topology.
\\
For all $n\in\N$ and $t\in[0;1]$, let us define $\varphi_{n,t}=id+t(\varphi_{n}-id)$.

For all $n\in\N$, we will denote by $H_{n}$ the diffeomorphism of $[0;1]$ defined by $H_{n}=\varphi_{n}\circ\ldots\circ\varphi_{0}$, and, if $T\in\R$, we define $h_{T}$,
diffeomorphism of $[0;1]$,
by : $h_{T}=\varphi_{n+1,t}\circ H_{n}$, where $t\in[0;1[$ represents the fractional part 
of $T$, and $n\in\N$ its integer part.
\\
Since 
$\varphi_{n,t}H_{n-1}\limite{t \to 1}{} \varphi_{n}H_{n-1}=H_{n}=id\circ H_{n}=\varphi_{n,0}H_{n}$, and 
because of the continuity of the path $(\varphi_{n,t})_{t\in[0;1]}$ for all $n\in\N$, 
the path $(h_{T})_{T\in\R}$ is continuous. 
By hypothesis, if $g\in G$, then $H_{n-1}gH_{n-1}^{-1}$ converges to $id$ when $n$ tends to 
$+\infty$ ; one knows also that $\Vert\varphi_{n,t}-id\Vert_{1}\leqslant\Vert\varphi_{n}-id\Vert_{1}$, where $\Vert\varphi_{n}-id\Vert_{1}$
tends to $0$ when $n$ tends to $+\infty$. Consequently, $\Vert\varphi_{n,t}-id\Vert_{1}$ tends uniformly to $0$ with respect to $t\in[0;1[$.
Thus $D(\varphi_{n,t}H_{n-1}gH_{n-1}^{-1}\varphi_{n,t}^{-1})$ tends uniformly to $1$ on $J$.
As $\varphi_{n,t}H_{n-1}gH_{n-1}^{-1}\varphi_{n,t}^{-1}$ has fixed points at the extremities of $J$, it follows that this diffeomorphism converges also to $id$ in $\mc{C}^{1}$-topology when $n$ tends to $+\infty$. Furthermore, if $g\in G$, we have the following equality :
$h_{T}gh_{T}^{-1}=\varphi_{n,t}H_{n-1}gH_{n-1}^{-1}\varphi_{n,t}^{-1}$, where $n$ still denotes the integer part of $T$ and $t$ its fractional part. As a consequence, $h_{T}gh_{T}^{-1}$ converges to $id$ when $T$ tends to $+\infty$.
\\
This completes the proof of the lemma.

\hfill$\Box$\medskip

\subsection{Groups which are isotopic to the identity by conjugacy are embeddable in a centralizer}

~~We show that a group which is isotopic to the identity by conjugacy is embeddable in a centralizer.

~~Let us consider a segment $J$
of $[0;+\infty[$, a group of diffeomorphisms of this segment $G$, and a
$\mc{C}^{1}$-continuous path $(h_{t})_{t\in[0;1[}$ of diffeomorphisms of $J$ such that $h_{0}=id$ and such that, 
for all $g\in G$, one has : $h_{t}gh_{t}^{-1}\limite{t\to\infty}{}id$.
\\
From Lemma \ref{l1.centralisateur et isotopie}, there exists a sequence $(\varphi_{n})_{n\in\N}$ of diffeomorphisms of $J$ converging to $id$ when $n$
tends to $+\infty$, and satisfying : 
$\varphi_{n}\ldots\varphi_{0}g\varphi_{0}^{-1}\ldots\varphi_{n}^{-1}
\limite{n\to\infty}{}id$, for all $g\in G$.
\\

~~Let now $f_{0}$ be a homothety with ratio $\alpha<1$ such that the iterates of $J$ by $f_{0}$ are pairwise disjoint.
One defines a contraction $f$ of $]0;+\infty[$ by the following : 
$\restriction{f}{J}=f_{0}\circ\varphi_{0}$ ;\\
$\restriction{f}{f_{0}(J)}=f_{0}\circ(f_{0}\varphi_{1}f_{0}^{-1})$ ; \ldots ;
$\restriction{f}{f_{0}^{n}(J)}=f_{0}\circ(f_{0}^{n}\varphi_{n}f_{0}^{-n})$
for all $n\geqslant 1$ ; $f=f_{0}$ elsewhere.
\\
Since $\varphi_{n}\limite{n\to\infty}{\mc{C}^{1}}id$, we know that $D(f_{0}^{n}\varphi_{n}f_{0}^{-n})(f_{0}^{n}(x)$ converges uniformly to $1$ with respect to $x\in J$.
Then
$\restriction{Df}{f_{0}^{n}(J)}$ tends to $\alpha$ when $n$ tends to $+\infty$, and thus $Df(x)\limite{x\to 0}{}\alpha$. Consequently, $f$ can be extended to $\mc{C}^{1}$-diffeomorphism of $[0;+\infty[$.
\\
Then we extend each $g_{J}\in G$ to a diffeomorphism $g$ of $]0;+\infty[$ in such a way that it commutes with $f$, i.e. as coinciding with
$g=f^{n}g_{J}f^{-n}$ on $f_{0}^{n}(J)$, where $n\in\Z$, and with $id$ elsewhere. 
For all $n\in \N$, we have then : 
$\restriction{g}{f_{0}^{n}(J)}=
f_{0}^{n}\varphi_{n}\ldots\varphi_{0}g_{J}\varphi_{0}^{-1}\ldots
\varphi_{n}^{-1}f_{0}^{-n}$, where $\varphi_{n}\ldots\varphi_{0}g_{J}\varphi_{0}^{-1}\ldots
\varphi_{n}^{-1}$ tends in $\mc{C}^{1}$-topology to the identity.
So, as $f_{0}$ is a homothety, the derivative of $\restriction{g}{f_{0}^{n}(J)}$ tends also to $1$ when $n$ tends to $+\infty$, and thus : $Dg(x)\limite{x\to 0}{}1$.
It follows that $g$ extends in a differentiable way at $0$.

\subsection{Groups which are embeddable in a centralizer are isotopic by conjugacy to the identity}

~~We show now the necessary condition of Theorem~\ref{t.1} : given a group $G$ of 
diffeomorphisms of some segment $J=[a;b]\subset[0;+\infty[$, one assumes that there exists a contraction $f$ of $[0;+\infty[$ in whose centralizer $G$
is embeddable.

~~First note that, by replacing if necessary $f$ by $f^{k}$, with $k$ some great enough integer, one can assume that $J$ is included in a fundamental domain $[f(x_{0});x_{0}]$ of $f$.

~~Each element $g_{J}$ of $G$ can by hypothesis extend in such a way that it commutes with $f$.
Denoting by $g$ this continuation, for all $n\in\N$ we have $\restriction{g}{f^{n}(J)}=f^{n}g_{J}f^{-n}$.
\\
If $n\in\N$, we denote by $A_{n}$ the increasing affine map from $f^{n}(J)$ into $J$ and 
by $g_{n}$ the $\mc{C}^{1}$-diffeomorphism of $J$ defined by : $g_{n}=A_{n}\restriction{g}{f^{n}(J)}A_{n}^{-1}$.
\\

\underline{Claim :} One has : $\Vert g_{n}-id\Vert_{1}\limite{n \to \infty}{}0$.
\\
Indeed, the orbits of fixed points of $g_{J}$ which are the extremities $a$ and $b$ of $J$ by $g$ accumulate in $0$ and are constituted with fixed points of $g$ ; by continuity of the derivative of $g$ at $0$, the derivative of $g$ tends to $Dg(0)$ at $0$ ; moreover $Dg(0)=\lim_{n \to\infty}\dfrac{g(f^{n}(a))}{f^{n}(a)}=1$.
Thus, if $n$ is a sufficiently greatwhole number, $\restriction{g}{f^{n}(J)}$ has derivative near  $1$ between two fixed points, and as a consequence is $\mc{C}^{1}$-close to $id$. By conjugating this diffeomorphism by the affine map $A_{n}$, one does not modify its derivative, which enables us to conclude that $g_{n}$ is also $\mc{C}^{1}$-close to the identity if $n$ is great enough.

\hfill$\Box$\medskip

From the expression of $\restriction{g}{f^{n}(J)}$ given above, one has : 
\[ g_{n}= A_{n}f^{n}g_{J}f^{-n}A_{n}^{-1}\]

Defining $\varphi_{0}=\restriction{id}{J}$ ; $\varphi_{1}=A_{1}\restriction{f}{J}$ ; $\varphi_{n}=(A_{n}\restriction{f^{n}}{J})(A_{n-1}\restriction{f^{n-1}}{f^{n-1}(J)})^{-1}$ for $n\geqslant 2$, one has : 
$g_{n}=\varphi_{n}\ldots\varphi_{1}\varphi_{0}g_{0}\varphi_{0}^{-1}\varphi_{1}^{-1}\ldots
\varphi_{n}^{-1}$. From the previous claim and Lemma \ref{l1.centralisateur et isotopie}, if we prove the following result, we will have proved that there exists an isotopy by conjugacy from $g_{J}$ to $id$ and by this way complete the proof of Theorem \ref{t.1} :
\\

\underline{Claim :} $\varphi_{n}\limite{n \to \infty}{}id$ with respect to the  $\mc{C}^{1}$-topology.
\\
Proof of the claim : By definition, one has : $\varphi_{n}=A_{n}\restriction{f}{f^{n-1}(J)}A_{n-1}^{-1}$.
Furthermore, as $Df$ is continuous at $0$ and $Df(0)\neq 0$, one has : $\sup\limits_{x,y\in f^{n-1}(J)}\dfrac{Df(x)}{Df(y)}\limite{n\to\infty}{}1$.
Consequently, $\dfrac{D\varphi_n(x)}{D\varphi_n(y)}$ tends uniformly to $1$ with respect to $x,y\in J$ when $n$ tends to $+\infty$. Since $D\varphi_n$ is equal to $1$ at at least one point of $J$, it follows that $D\varphi$ tends uniformly to $1$ on $J$. The claim follows.

\hfill$\Box$\medskip

\section{Isotopies by conjugacy}
 
\subsection{Statement of the result}
~~In this section, we consider the set 
\\
 $D_{\alpha,\beta}=\left\lbrace f\in Diff^{1}(\left[ 0;1\right] ) ; Fix(f)=\left\lbrace 0;1\right\rbrace , f\geqslant id , Df(0)=\alpha , Df(1)=\beta \right\rbrace$, well-defined for all real numbers $\alpha\geqslant 1$ and $0<\beta\leqslant 1$, and we show that the classes of conjugacy of such diffeomorphisms are dense in this set for the  $\mc{C}^{1}$-topology.
\\

~~Let now $f$ and $g$ be two diffeomorphisms of $D_{\alpha,\beta}$. We will show that in each $\mc{C}^{1}$-neighbourhood of $g$, there is a conjugate of $f$. Thus, we will have showed :

\begin{theo}\label{t.densité C1} For all diffeomorphism $f\in \D $, the differentiable class of conjugacy of $f$ is dense in $D_{\alpha,\beta}$.
\end{theo}

Let $f,g$ be two elements of $\D$, and $\mc{U}$ be a $\mc{C}^{1}$-neighbourhood of $g$ in $\D$.
In order to find a conjugate of $f$ in $\mc{U}$, we will perturb $g$ by sufficiently small diffeomorphisms so that we don't go out $\mc{U}$, till we find a conjugate of $f$.

As announced in the introduction, we actually show the stronger result stated in Theorem \ref{t.4}, which we will prove from the following theorem :

\begin{theo}\label{t.transfochemin} Let $f$ be a diffeomorphism of $D_{\alpha,\beta}$ without fixed point in $]0,1[$ and $(\varepsilon_t)_t\in[0;1[$ be a continuous path of strictly positive numbers. Let $(f_t)_{t\in[0,1]}$ be a $\mc{C}^1$-continuous path of diffeomorphisms of $D_{\alpha,\beta}$ such that :
\begin{itemize}
\item $f_0=f$ ;
\item $f_t$ has no fixed point in $]0,1[$ for all $t<1$.
\end{itemize}
Then there exists a $\mc{C}^1$-continuous path $(h_t)_{t\in[0,1[}$ of diffeomorphisms of $[0,1]$, all of them coinciding with $id$ on a neighbourhood of $0$ and of $1$, such that $h_0=id$ and $\Vert h_t f h_t^{-1}-f_t\Vert_1<\varepsilon_t$ for all $t\in[0;1[$.
\end{theo}

\underline{Remarks :} --Since $f_t$ has no fixed point in $]0,1[$ for all $t<1$, the diffeomorphism $f_{1}$ can have in $]0;1[$ only fixed points with derivative equal to $1$. 
\\
--By choosing the continuous path $(\varepsilon_t)$ as converging to $0$, under the hypotheses of this theorem, one obtains a continuous path $(h_t)_{t\in[0,1[}$ such that
$\lim\limits_{t\to 1} h_tfh_t^{-1} =\lim\limits_{t\to 1} f_t$.
\\

\underline{Proof of Theorem \ref{t.4} by Theorem \ref{t.transfochemin} :}
Given two diffeomorphisms $f$ and $g$ in $D_{\alpha,\beta}$, one can always easily exhibit a $\mc{C}^{1}$-continuous path $(f_t)_{t\in[0;1[}$ of diffeomorphisms of $D_{\alpha,\beta}$, such that $f_0=f$ et $f_t\limite{t\to 1}{}g$.
For example, let $f_t$ be defined by $f_t=(1-t)f+tg$.
If $(\varepsilon_t)_{t\in[0;1[}$is a continuous path of strictly positive numbers converging to $0$, Theorem \ref{t.transfochemin} provides a continuous path $(h_t)_{t\in[0,1[}$ of diffeomorphisms of $[0,1]$ such that 
$h_0fh_0^{-1}=f$ ; $h_tfh_t^{-1}\limite{t\to 1}{}g$ and such that $(h_t)_{t\in[0;1[}$ is continuous. Then, 
$(h_tfh_t^{-1})_{t\in[0;1[}$ is an isotopy by conjugacy from $f$ to $g$.

\hfill$\Box$\medskip


\subsection{Two "gluing lemmas"}\label{s.recollement}

~~From now on and in all the following, $\Phi$ will denote the $\mc{C}^{\infty}$-diffeomorphism of $[0;1]$ which we define now as being 
decreasing, being equal to $1$ on $[0;\frac{1}{2}]$ and equal to $0$ on a neighbourhood of $1$. Its derivative, continuous on the compact set $[0;1]$, is consequently bounded. We define $M_{\Phi}=\max\limits_{x\in[0;1]}\vert D\Phi\vert>1$.

\paragraph{gluing at the extremities}

~~Given two $\mc{C}^{1}$-diffeomorphisms of the interval, $f$ and $g$, having same derivatives at $0$ and at $1$, and a strictly positive real number $\varepsilon$, Corollary \ref{c1.densité C2} enables us to obtain a new $\mc{C}^{1}$-diffeomorphism which coincides with $f$ on a neighourhood of $0$ and on a neighourhood of $1$, and which is $\varepsilon$-close to $g$ for the $\mc{C}^{1}$-topology on the whole interval $[0;1]$.

\begin{defi}\label{d.densité C2}
Let $g$ be an element of $\mc{D}iff^{1}([0;1])$. For all strictly positive real number $\varepsilon$, we will denote by $\mc{U}_{\varepsilon}^{g}$ the set of maps 
$f\in\mc{D}iff^1([0;1])$ such that there exist two real numbers $a$ and $b$ of $]0;1[$ satisfying :
\begin{itemize}
\item [$\bullet$] $\Vert \restriction{f}{[0;a]}-\restriction{g}{[0;a]}\Vert_{1}<\varepsilon$
\item [$\bullet$] $\Vert \restriction{f}{[b;1]}-\restriction{g}{[b;1]}\Vert_{1}<\varepsilon$.
\end{itemize}
If $a$ and $b$ are two real numbers of $]0;1[$, we will denote by $\mc{U}_{\varepsilon,a,b}^{g}$ the set of maps 
$f\in\mc{D}iff^1([0;1])$ satisfying the two properties above.
\end{defi}

\begin{lem}\label{l2.densité C2} Let $\varepsilon>0$ be a real number.
\\
Then there exists $\tilde{\varepsilon}>0$ such that : 
for all $g\in \mc{D}iff^{1}([0;1])$, for all $a,b\in]0;1[$, 
if $f\in\mc{U}_{\tilde{\varepsilon},a,b}^{g}$,
then there exists $f_{0}\in \mc{D}iff^{1}([0;1])$ such that :
\begin{enumerate}
\item $\Vert f_{0}-g\Vert_{1}<\varepsilon$
\item $\restriction{f_{0}}{[0;\frac{a}{2}]\cup[\frac{1+b}{2};1]}
=\restriction{f}{[0;\frac{a}{2}]\cup[\frac{1+b}{2};1]}$
\item $\restriction{f_{0}}{[a;b]}=\restriction{g}{[a;b]}$
\end{enumerate}

\end{lem}

\underline{Proof :} 

~~Let $\varepsilon>0$ be a real number, and $0<\tilde{\varepsilon}<\dfrac{\varepsilon}{2M_{\Phi}}$.

Let $g$ be an element of $\mc{D}iff^{1}([0;1])$ ; $a,b\in]0;1[$ and $f\in \mc{U}_{\tilde{\varepsilon},a,b}^{g}$. If $a\geqslant b$, it suffices to let $f_{0}$ be $f$ to conclude. Thus, let us consider now the case where $a<b$.

~~We define now the map $\Phi_{0}$ of $[0;1]$, which is obtained from $\Phi$ as follows : 
\begin{itemize}
\item [$\bullet$] if $x\in [0;a]$, then $\Phi_{0}(x)=\Phi(\frac{x}{a})$. 
So $\Phi_{0}$ is equal to $1$ on $[0;\frac{a}{2}]$ and to $0$ on a neighbourhood of $a$ in $[0;a] $ ; 

\item[$\bullet$] $\Phi_{0}(x)=0$ if $x\in 
[a;b] $ ;

\item[$\bullet$] $\Phi_{0}(x)=\Phi(\frac{-x+1}{1-b})$ if $x\in[b;1]$, so $\Phi_{0}$ is equal to $0$ on a neighbourhood of $b$ in $[b;1]$ and to $1$ on 
$[\frac{1+b}{2};1]$.
\end{itemize}

Let us now define $f_{0}=\Phi_{0}f+(1-\Phi_{0})g$, and
let us notice that the hypothesis
"$\Vert \restriction{f}{[0;a]}-\restriction{g}{[0;a]}\Vert_{1}<\tilde{\varepsilon}$" implies, by integration of the inequality 
$\vert Df-Dg\vert<\tilde{\varepsilon}$ on $[0;x]$, that : $\vert f(x)-g(x)\vert<\tilde{\varepsilon}x$ for all $x\in[0;a]$. By the same method, one obtains the inequality $\vert f(x)-g(x)\vert<\tilde{\varepsilon}(1-x)$, satisfied by each $x\in[b;1]$.

As the conclusions 2. et 3. of the lemma follows immediately from the construction of $f_{0}$,
one has now only to show that the map $f_{0}$ satisfies : 
$\vert f_{0}-g\vert <\varepsilon$ and $ \vert Df_{0}-Dg\vert <\varepsilon$, which can be stated by simple calculation. 

\hfill$\Box$\medskip

~~From this lemma, we deduce the following corollary :

\begin{coro}\label{c1.densité C2}
Let $\varepsilon>0$ be a real number, and $f,g\in D_{\alpha,\beta}$. Then there exist $a,b\in]0;1[$ and $f_{0}\in D_{\alpha,\beta}$ such that :
\begin{enumerate}
\item $\Vert f_{0}-g\Vert_{1}<\varepsilon$
\item $\restriction{f_{0}}{[0;\frac{a}{2}]\cup[\frac{1+b}{2};1]}
=\restriction{f}{[0;\frac{a}{2}]\cup[\frac{1+b}{2};1]}$
\item $\restriction{f_{0}}{[a;b]}=\restriction{g}{[a;b]}$.
\end{enumerate}
\end{coro}

\underline{Proof :} Indeed, considering $\tilde{\varepsilon}>0$ as in Lemma \ref{l2.densité C2}, since $f$ and $g$ have same derivatives at $0$ and at $1$, there exist $a,b\in]0;1[$ such that $f\in\mc{U}_{\tilde{\varepsilon},a,b}^{g}$.
It suffices then to apply Lemma \ref{l2.densité C2} to conclude.

\hfill$\Box$\medskip

~~This result is the one that will be usefull to prove the only density of classes of conjugacy of diffeomorphisms of $D_{\alpha,\beta}$ in that set ; however, in order to obtain an isotopy by conjugacy from a diffeomorphism of $D_{\alpha,\beta}$ to another, a  parameter version of it will be needed :

\begin{lem}\label{l2.densité C2 bis}
Let us consider $f\in D_{\alpha,\beta}$ and a $\mc{C}^{1}$-continuous path $(f_{t})_{t\in[0;1[}$ of diffeomorphisms of $D_{\alpha,\beta}$.
\\
Then, for all continuous path of strictly positive real numbers $(\varepsilon_{t})_{t\in[0;1[}$, there exist two continuous paths $(a_{t})_{t\in[0;1[}$ and $(b_{t})_{t\in[0;1[}$ of real numbers of $]0;1[$ ; there exists a new $\mc{C}^{1}$-continuous path $(\tilde{f}_{t})_{t\in[0;1[}$ of diffeomorphisms of $D_{\alpha,\beta}$ such that :
\begin{itemize} 
 \item for all $t<1$, $\restriction{\tilde{f}_{t}}{[0;\frac{a_{t}}{2}]\cup[\frac{b_{t}+1}{2};1]}=
\restriction{{f}}{[0;\frac{a_{t}}{2}]\cup[\frac{b_{t}+1}{2};1]}$,
\item for all $t<1$, $\Vert\tilde{f}_{t}-f_{t}\Vert_{1}<\varepsilon_{t}$.
\end{itemize}
Moreover, if $f_0=f$, then we can also demand that $\tilde{f}_0=f$.
\end{lem}

Let us give first the following result :

\begin{lem}\label{l3.densité C2 bis}
Let $(\varepsilon_{t})_{t\in[0;1[}$ be a continuous path of strictly positive real numbers.
Consider $f\in D_{\alpha,\beta}$ and a $\mc{C}^{1}$-continuous path $(f_{t})_{t\in[0;1[}$ of diffeomorphisms of $D_{\alpha,\beta}$.
\\
Then there exist two continuous paths $(a_{t})_{t\in[0;1[}$ and $(b_{t})_{t\in[0;1[}$ of real numbers of $]0;1[$ such that, for all $t<1$, 
one has $\Vert \restriction{f_{t}}{[0;a_{t}]\cup[b_{t};1]}-\restriction{f}{[0;a_{t}]\cup[b_{t};1]}\Vert_{1}<\varepsilon_{t}$.
\end{lem}

\underline{Proof of Lemma \ref{l3.densité C2 bis} :}
Given an increasing sequence $(t_n)_{n\in\N}$ of $[0;1[$ converging to $1$, one uses the continuity of the path $(f_t)_{t\in[0;1[}$ on each compact set $[t_n;t_{n+1}]$ so as to obtain a partition of this interval in finitely many subintervals such that, for all $t$ in one of these subintervals $I$, an only real number $0<a_I<1$ satisfies : $\Vert \restriction{f_t}{[0;a_I]}-\restriction{f}{[0;a_I]}\Vert_1<\min\limits_{t\in[t_{n};t_{n+1}]}\varepsilon_t$. All what is needed now is finding a continuous map $\varphi$ such that, for all $t$ in a subinterval $I$, $a_t=\varphi(t)$ is smaller than $a_I$, which can easily be done.
\\
We conclude by the same token as regards the path $(b_t)$.

\hfill$\Box$\medskip

\underline{Proof of Lemma \ref{l2.densité C2 bis} :}
Lemma \ref{l3.densité C2 bis} gives us two continuous paths $(a_t)$ and $(b_t)$ such that for all $t$, 
$\Vert \restriction{f_t}{[0;a_t]\cup[b_t;1]}-\restriction{f}{[0;a_t]\cup[b_t;1]}\Vert_1
<\inf(\frac{\varepsilon_t}{2};\frac{\varepsilon_t}{4M_\Phi})$.
Considering the construction done in Lemma \ref{l2.densité C2}, one obtains the wanted path $(\tilde{f}_t)$ : one checks that the continuity of this path follows from the continuity of the paths $(a_t)$ and $(b_t)$.
 
\hfill$\Box$\medskip

To show Theorem \ref{t.transfochemin}, it is sufficient now to show the following theorem :

\begin{theo}\label{t1.transfochemin} Let $f$ be a diffeomorphism of $D_{\alpha,\beta}$ without fixed point in $]0,1[$ and $(\varepsilon_t)_t\in[0;1[$ be a continuous path of strictly positive real numbers. Let $(f_t)_{t\in[0,1]}$ be a $\mc{C}^1$-continuous path of diffeomorphisms of $D_{\alpha,\beta}$ such that :
\begin{itemize}
\item $f_0=f$ ;
\item $f_t$ has no fixed point in $]0,1[$ for all $t<1$ ;
\item For all $t<1$, there exist $a_t>0$ and $b_t<1$, continuously depending on $t$, such that $f_t$ coincides with $f$ on $[0,a_t]\cup[b_t,1]$.
\end{itemize}
Then there exists a continuous path $(h_t)_{t\in[0,1[}$ of diffeomorphisms of $[0,1]$, all of them coinciding with $id$ on a neighbourhood of $0$ and of $1$, such that $h_0=id$ and $\Vert h_t f h_t^-1-f_t\Vert_1<\varepsilon_t$ for all $t\in[0;1[$.
\end{theo}

\underline{Proof of Theorem \ref{t.transfochemin} by Theorem \ref{t1.transfochemin} :}
From a continuous path $(f_t)_{t\in[0;1[}$ of diffeomorphisms of $D_{\alpha,\beta}$ and thanks to Lemma \ref{l2.densité C2 bis}, one can exhibit two continuous paths $(a_t)_{t\in[0;1[}$ and $(b_t)_{t\in[0;1[}$ and a continuous path of diffeomorphisms $(\tilde{f}_t)_{t\in[0;1[}$ such that :
\begin{itemize} 
\item $\tilde{f}_0=f_0$ ;
\item for all $t<1$, $\restriction{\tilde{f}_{t}}{[0;\frac{a_{t}}{2}]\cup[\frac{b_{t}+1}{2};1]}=
\restriction{{f}}{[0;\frac{a_{t}}{2}]\cup[\frac{b_{t}+1}{2};1]}$ ;
\item for all $t<1$, $\Vert\tilde{f}_{t}-f_{t}\Vert_{1}<\dfrac{\varepsilon_{t}}{2}$.
\end{itemize}
One applies then Theorem \ref{t1.transfochemin} to obtain a $\mc{C}^{1}$-continuous path $(h_{t})_{t\in[0;1[}$ of diffeomorphisms of $[0;1]$, with $h_0=id$, such that $h_t$ coincides with $id$ on a neighbourhood of $0$ and of $1$, and such that $\Vert h_t f h_t^-1-\tilde{f}_t\Vert_1<\dfrac{\varepsilon_t}{2}$ for all $t\in[0;1[$.
Then we have : $\Vert h_t f h_t^-1-f_t\Vert_1<\varepsilon_t$ for all $t\in[0;1[$.
\hfill$\Box$\medskip

~~Let us also notice that this corollary of Lemma \ref{l2.densité C2} which will be usefull in section \ref{s.point fixe dégénéré}. If $f$ and $g$ are two $\mc{C}^{1}$-diffeomorphisms  of $[0;1]$ coinciding at one point, this corollary enables us to exhibit a new $\mc{C}^{1}$-diffeomorphism coinciding with $f$ in a neighbourhood of this point, and arbitrarily close to $g$.

\begin{coro}\label{c.densité C1}
Let $\varepsilon>0$ be a real number. Then there exists $\tilde{\varepsilon}>0$ such that : for all $g\in\mc{D}iff^1([0;1])$, for all point $x_0$ of $]0;1[$, for all $\eta\in]0;\min(\frac{x_0}{2},\frac{1-x_0}{2})[$, if $f$ is a diffeomorphism of $[0;1]$ such that
\begin{enumerate}
\item $f(x_0)=g(x_0)$
\item $\Vert\restriction{f}{[x_0-2\eta;x_0+2\eta]}-\restriction{g}{[x_0-2\eta;x_0+2\eta]}\Vert_1
<\tilde{\varepsilon}$,
\end{enumerate}
then there exists $g_0\in\mc{D}iff^1([0;1])$ such that :
\begin{enumerate}
\item $\restriction{g_0}{[x_0-\eta;x_0+\eta]}=\restriction{f}{[x_0-\eta;x_0+\eta]}$
\item $\restriction{g_0}{[0;x_0-2\eta]\cup[x_0+2\eta;1]}=\restriction{g}{[0;x_0-2\eta]\cup[x_0+2\eta;1]}$
\item $\Vert g_0-g\Vert_1<\varepsilon$.
\end{enumerate}
\end{coro}

Here follows the corresponding parameter version of this corollary :

\begin{coro}\label{c.densité C1 bis} \label{c.transfochemin} Let $(\varepsilon_t)_{t\in[0,1[}$ be a continuous path of strictly positive real numbers. 
Then there exists a continuous path $(\tilde{\varepsilon}_t)_{t\in[0,1[}$ such that : for all $\mc{C}^1$-continuous path $(f_t)_{t\in[0;1[}$ of diffeomorphisms of $\mc{D}iff^1([0;1])$, for all continuous path of points $(x_t)_{t\in[0,1[}$ of $]0;1[$, for all $\eta_n\in]0;\min(\frac{x_t}{2},\frac{1-x_t}{2})[$, if $f$ is a diffeomorphism of $[0;1]$ such that
\begin{enumerate}
\item $f(x_t)=f_t(x_t)$ for all $t<1$
\item $\Vert\restriction{f}{[x_t-2\eta_t;x_t+2\eta_t]}-\restriction{f_t}{[x_t-2\eta_t;x_t+2\eta_t]}\Vert_1
<\tilde{\varepsilon_t}$,
\end{enumerate}
then there exists a $\mc{C}^1$-continuous path $(g_t)_{t\in[0,1[}$ of $\mc{D}iff^1([0;1])$, such that :
\begin{enumerate}
\item $\restriction{g_t}{[x_t-\eta_t;x_t+\eta_t]}=\restriction{f}{[x_t-\eta_t;x_t+\eta_t]}$
\item $\restriction{g_t}{[0;x_t-2\eta_t]\cup[x_t+2\eta_t;1]}=\restriction{f_t}{[0;x_t-2\eta_t]\cup[x_t+2\eta_t;1]}$
\item $\Vert g_t-f_t\Vert_1<\varepsilon_t$.
\end{enumerate}

\end{coro}
Prooves of Corollaries~\ref{c.densité C1}~and~\ref{c.densité C1 bis} are absolutely similar to the one of Lemma \ref{l3.densité C2 bis}, thus we omit it.

\paragraph{"partial gluing" near to an extremity}

\begin{defi}
If $g\in\mc{D}iff^{1}([0;1])$, $\varepsilon$ is a strictly positive real number and
$a\in]0;1[$, 
we willl denote by $\mc{U}_{\varepsilon,a}^{g}$ the set of maps $f\in\mc{H}omeo([0;1])$ such that : 
\begin{itemize}
\item[$\bullet$] $\restriction{f}{[a;1]}\in\mc{D}iff^{1}([a;1],[f(a);1])$
\item[$\bullet$] $\Vert \restriction{f}{[a;1]}-\restriction{g}{[a;1]}\Vert_{1}<\varepsilon$.
\end{itemize}
\end{defi}

~~We will see in the proof of the following lemma how, from a $\mc{C}^{1}$-diffeomorphism $g$ of $[0;1]$ and from a homeomorphism $f\in\mc{U}_{\varepsilon,a}^{g}$, one can exhibit a map coinciding with $f$ in a neighbourhood of the extremities and which is arbitrarily close to $g$.

\begin{lem}\label{l2.densité C1}
Let $\varepsilon>0$ be a real number.
Then there exists $\tilde{\varepsilon}>0$ such that : 
for all $g\in\mc{D}iff^{1}([0;1])$, for all $a\in]0;1[$ 
and $1>b>\dfrac{a+1}{2}$, for all $f\in\mc{U}_{\tilde{\varepsilon},a}^{g}$, 
there exists $g_{0}\in\mc{D}iff^{1}([0;1])$ satisfying :
\begin{itemize}
\item [$\bullet$] $\restriction{g_{0}}{[\frac{a+1}{2};b]}=\restriction{f}{[\frac{a+1}{2};b]}$
\item [$\bullet$] $\restriction{g_{0}}{[0;a]\cup[\frac{b+1}{2};1]}
=\restriction{g}{[0;a]\cup[\frac{b+1}{2};1]}$
\item [$\bullet$] $\Vert g_{0}-g\Vert_{1}<\varepsilon$.
\end{itemize}
\end{lem}

\underline{Proof of Lemma \ref{l2.densité C1} :}
Let $\varepsilon>0$ be a real number, 
and $\tilde{\varepsilon}<\min\left(\dfrac{\varepsilon}{2};\dfrac{\varepsilon}{4M_{\Phi}} \right)$.
Let us consider $g\in\mc{D}iff^{1}([0;1])$, $a\in]0;1[$, $f\in\mc{U}_{\tilde{\varepsilon},a}^{g}$, and $\Phi_{0}$,
$\mc{C}^{\infty}$-map of $[0;1]$, constructed from $\Phi$ (see Section \ref{s.recollement})
by the following way :
\begin{itemize}
\item [$\bullet$] on $[0;a]$ and on $[\frac{b+1}{2};1]$, $\Phi_{0}$ is constant, equal to $1$ ;
\item[$\bullet$] if $x\in[a;\frac{a+1}{2}]$, 
then $\Phi_{0}(x)=\Phi\left(\dfrac{x}{\frac{1-a}{2}}\right)$ ;
\item[$\bullet$] on $[\frac{a+1}{2};b]$, $\Phi_{0}$ is constant, equal to $0$ ;
\item[$\bullet$] if $x\in[b;\frac{b+1}{2}]$, 
then $\Phi_{0}(x)=\Phi\left(\dfrac{1-x}{\frac{1-b}{2}}\right)$.
\end{itemize}

Then, one defines $g_{0}=\Phi_{0}\,g+(1-\Phi_{0})f$, and one checks that
$\Vert g-g_{0}\Vert_{1}$ is bounded by $\varepsilon$. The other conclusions are satisfied by construction :
one can check it by simple calculation thanks to this remark : the integration of the inequality $\vert Df-Dg\vert<\tilde{\varepsilon}$ gives : 
\\
$\vert f(x)-g(x)\vert<\tilde{\varepsilon}(1-x)$, for all $x\in[a;1]$.

\hfill$\Box$\medskip

\underline{Remark :}
One can notice that, in that Lemma \ref{l2.densité C1} as in Lemma \ref{l2.densité C2}, $\tilde{\varepsilon}$ dépends only on $\varepsilon$. That justifies the following definition :

\begin{defi} Given a real number $\varepsilon>0$, we will denote by \emph{$\textrm{Marg}(\varepsilon)$} 
the set of strictly positive real numbers satisfying the conditions satisfied by $\tilde{\varepsilon}$ in Lemmas \ref{l2.densité C2} and \ref{l2.densité C1} and in Corollary \ref{c.densité C1}. In other words : 
$\textrm{Marg}(\varepsilon)=]0;\min(\frac{\varepsilon}{2},\frac{\varepsilon}{4M_{\Phi}})[$
\end{defi}

As we did for the first "gluing" lemma, one states here again a parameter version of this lemma :

\begin{lem}\label{l2.densité C1 bis}
Let $(\varepsilon_t)_{t\in[0;1[}$ be a continuous path of strictly positive real numbers, $f\in D_{\alpha,\beta}$ and $(f_t)_{t\in[0;1[}$ a continuous path of diffeomorphisms of $D_{\alpha,\beta}$.
\\
Lemma \ref{l2.densité C2 bis} ensures then the existence of a continuous path $(b_t)_{t\in[0;1[}$ such that $\Vert \restriction{f_t}{[b_t;1]}-\restriction{f}{[b_t;1]}\Vert_{1}<\textrm{Marg}(\varepsilon_t)$.
Then, for all paths $(c_t)_{t\in[0;1[}$ and $(d_t)_{t\in[0;1[}$ such that $b_t<c_t<1$ and $\dfrac{c_t+1}{2}<d_t<1$, there exists a $\mc{C}^{1}$-continuous path $(\tilde{f}_t)_{t\in[0;1[}$ of diffeomorphisms of $D_{\alpha,\beta}$ satisfying :
\begin{enumerate}
\item $\restriction{\tilde{f}_t}{[0;c_t]\cup[\frac{d_t+1}{2};1]}=\restriction{f_t}{[0;c_t]\cup[\frac{d_t+1}{2};1]}$,
\item $\restriction{\tilde{f}_t}{[\frac{c_t+1}{2};d_t]}=\restriction{f}{[\frac{c_t+1}{2};d_t]}$,
\item $\Vert\tilde{f}_t-f_t\Vert_{1}<\varepsilon_t$.
\end{enumerate}
\end{lem}

Once again, the proof is analogous to the one of Lemma~\ref{l3.densité C2 bis}, and thus we omit it.

\hfill$\Box$\medskip

\subsection{Unitary conjugacy from $f$ to $g_0$}

If $f\in D_{\alpha,\beta}$, we introduce the following notations : 
\\
$D_{f}=\lbrace g\in D_{\alpha,\beta} ; \exists a,b,
\restriction{f}{\left[ 0;a\right]} =\restriction{g}{\left[ 0;a\right]} \mbox{ and } \restriction{f}{\left[ b;1\right]}=\restriction{g}{\left[ b;1\right]}\rbrace$ ;
\\
and, if $a,b\in \left] 0;1\right[$ :
\\ 
$D_{f,a,b}=\lbrace g\in D_{\alpha,\beta} ; \restriction{f}{\left[ 0;a\right]} =\restriction{g}{\left[ 0;a\right]} \mbox{ and } \restriction{f}{\left[ b;1\right]}=\restriction{g}{\left[ b;1\right]}\rbrace$.

~~If $f$ is a diffeomorphism of $D_{\alpha,\beta}$, $a$ and $b$ are two real numbers such that $1>a,b>0$ and $g$ belongs to $D_{f,a,b}$ ; if $x_{0}<a$ is such that $f(x_{0})<a $, 
then : $\restriction{f}{\left[ 0 ; f(x_{0})\right] }=\restriction{g}{\left[ 0 ; f(x_{0})\right] }$, and then 
$h$ defined on $\left[ 0 ; f(x_{0})\right]$ by $h=id $ conjugates $f$ to $g$ on $\left[ 0 ; x_{0}\right] $. 
Remark below states that this diffeomorphism extends uniquely on $\left[0;1\right[$ in such a way that it conjugates $f$ to $g$.
\\

\underline{Remark :} For all $f\in D_{\alpha,\beta}$ and $g\in D_{f}$, there exists an unique diffeomorphism $h$ of $[0;1[$ which conjugates $f$to $g$ and coinciding with $id$ on a neighbourhood of $0$.

\underline{Proof :}
If $h$ conjugates $f$ to $g$, the relation $hfh^{-1}=g$ gives, for all integer $n$, $\restriction{h_{n}}{\left[ f^{n}(x_{0}) ; f^{n+1}(x_{0})\right]}=g^{n}h_{0}f^{-n}=g^{n}f^{-n}$. One checks that this formula defines a $\mc{C}^1$-diffeomorphism of $[0;1[$.

\hfill$\Box$\medskip

\begin{defi} The homeomorphism of $[0;1]$ coinciding with the diffeomorphism $h$ mentioned in the remark above on $[0;1[$ and with value $1$ at $1$ will be called in the following \emph{unitary conjugacy from $g$ to $f$}, and denoted by $h_{g}$. It is the only conjugacy from $g$ to $f$ which coincides with $id$ on a neighbourhood of $0$.
\end{defi}

\begin{lem}\label{l1.conjugante}
Let $f$ be a $\mc{C}^1$-diffeomorphism of $[0;1]$ ; $(a_t)_{t\in[0;1[}$ and $(b_t)_{t\in[0;1[}$ be two continuous paths of real numbers of $[0;1]$ and $(f_t)_{t\in[0;1[}$ be a $\mc{C}^1$-continuous path of diffeomorphisms of $[0;1]$ such that, for all $t\in[0;1[$, we have $f_t\in D_{f,a_t,b_t}$.
\\
Then, for all compact interval $I$ of $]0;1[$, the path $(\restriction{h_{f_t}}{I})_{t\in[0;1[}$ is continuous with respect to the $\mc{C}^1$-topology.
\end{lem}

\underline{Proof of Lemma \ref{l1.conjugante} :}
One denotes $I=[v_0;v_1]$. Since $v_1<1$, and by continuity of the path $(a_t)_{t}$, for all $t\in [0;1[$, there exists a neighbourhood $\mc{V}_t$ of $t$ and an integer $n_t$ such that, for all $t'\in\mc{V}_t$, $h_{f_{t'}}$ is defined by 
$h_{f_{t'}}=f_{t'}^{n_t}f^{-n_t}$, and that diffeomorphism continuously depends on $f_{t'}$ with respect to the $\mc{C}^{1}$-topology.
The continuity of $\restriction{h_{f_t}}{I}$ follows.
\hfill$\Box$\medskip

~~If the extension of $h_{g}$ to $[0;1]$ has $\mc{C}^{1}$-regularity, 
then $g$ is $\mc{C}^{1}$-conjugate to $f$ by $h_{g}$, and the conclusion of Theorem
\ref{t.densité C1} is satisfied. However, generally, $h_{g}$ is not differentiable at $1$ ; in order to show that the conjugacy class of $f$ is arbitrarily close to $g$, one decides thus to modify $g$ on $\mc{U}\cap D_{f}$ in such a way that $h_{g}$ coincides with $id$ on a neighbourhood of $1$.

\begin{prop}\label{p1.densité C1}
If $f\in D_{\alpha,\beta}$, $g\in D_{f}$, and $\mc{U}$ is a $\mc{C}^{1}$-neighbourhood of $g$, then there exists $g'\in D_{f}\cap\mc{U}$ such that
$h_{g'}=h'$ coincides with $id$ in a neighbourhood of $1$.
\end{prop}

~~The parameter version of this statement is the following : 

\begin{prop}\label{p1.densité C1 para}
Let us consider $f\in D_{\alpha,\beta}$ and $(a_t)_{t\in[0;1[}$ ; $(b_t)_{t\in[0;1[}$ two continuous paths of real numbers of $]0;1[$ ; $(\varepsilon_t)_{t\in[0;1[}$ a continuous path of strictly positive real numbers.
Let lastly $(f_t)_{t\in[0;1[}$ be a $\mc{C}^1$-continuous path of diffeomorphisms of $D_{\alpha,\beta}$, such that, for all $t\in[0;1[$, $f_t$ belongs to $D_{f,a_t,b_t}$.
\\
Then there exist two continuous paths $(\tilde{b}_t)_{t\in[0;1[}$ and $(\tilde{c}_t)_{t\in[0;1[}$ of real numbers of $]0;1[$, as well as a $\mc{C}^1$-continuous path $(\tilde{f}_t)_{t\in[0;1[}$ of diffeomorphisms of $D_{f,a_t,\tilde{b}_t}$ such that, for all $t\in[0;1[$, we have :
\begin{itemize}
\item $h_{\tilde{f}_t}=id$ sur $[\tilde{c}_t;1]$ ;
\item $\Vert \tilde{f}_t-f_t\Vert_{1}<\varepsilon_t$.
\end{itemize}
\end{prop}

In particular, this corollary follows : 

\begin{coro}\label{c2. densité C1}
Under the hypotheses of this proposition, the map $\begin{array}{ccc}
[0;1[ & \longrightarrow & \mc{D}iff^1([0;1])\\
t & \longmapsto & h_{\tilde{f}_t}
\end{array}$ is continuous with respect to the $\mc{C}^1$-topology.
\end{coro}

~~Indeed, we saw in Lemma \ref{l1.conjugante} that, for all compact interval $I$ of $]0;1[$, the diffeomorphism $\restriction{h_{\tilde{f}_t}}{I}$ continuously depends on $t$ with respect to the $\mc{C}^1$-topology. Thus it is sufficient to choose $I$ in such a way that 
$h_{\tilde{f}_t}$ coincides with $id$ on the complement of $I$ so as to ensure that
$\Vert \restriction{h_{\tilde{f}_{t'}}}{I}-\restriction{h_{\tilde{f}_t}}{I}\Vert_{1}$, as well as $\Vert \restriction{h_{\tilde{f}_{t'}}}{\complement{I}}
-\restriction{h_{\tilde{f}_t}}{\complement{I}}\Vert_{1}$ are arbitrarily small provided that $t$ and $t'$ are sufficiently close to each other.

~~From these two propositions, one can prove rather easily Theorems~\ref{t.densité C1} and \ref{t1.transfochemin}. Proposition \ref{p1.densité C1 para} will also be usefull in section \ref{s.isotopie identité} in the proof of Proposition \ref{p1.isotopie à id}.
\\

\underline{Proof of Theorem \ref{t.densité C1} by Proposition \ref{p1.densité C1} :}
Given two diffeomorphisms $f,g$ of $D_{\alpha,\beta}$, we prove that, arbitrairily $\mc{C}^1$-close to $g$, there exists a conjugate of $f$.
Let $\varepsilon$ be a strictly positive real number.
From Corollary \ref{c1.densité C2}, there exist $1>a,b>0$ and $g_{0}\in D_{f,a,b}$ such that $\Vert g_{0}-g\Vert_{1}<\dfrac{\varepsilon}{2}$.
From Proposition \ref{p1.densité C1}, there exists $\tilde{g}_0\in D_f$ such that 
$\Vert \tilde{g}_0-g_0\Vert_{1}<\dfrac{\varepsilon}{2}$ and $h_{\tilde{g}_0}=id$ in a neighbourhood of $1$.
Then $h_{\tilde{g}_0}$ is a $\mc{C}^1$-diffeomorphism of $[0;1]$ coinciding with the identity on a neighbourhood of $0$ and of $1$, conjugating $f$ to $\tilde{g}_0$ ; and $\Vert \tilde{g}_0-g\Vert_{1}<\varepsilon$.
\hfill$\Box$\medskip

\underline{Proof of Theorem \ref{t1.transfochemin} by Proposition \ref{p1.densité C1 para} :}
If $(f_t)_{t\in[0;1[}$ is a continuous path of diffeomorphisms of $D_{f}$, then the path $(\tilde{f}_t)_{t\in[0;1[}=(h_{\tilde{f}_t}fh_{\tilde{f}_t}^{-1})_{t\in[0;1[}$ given by Proposition \ref{p1.densité C1 para} is a $\mc{C}^{1}$-continuous path of conjuguates of $f$ such that the conjugacies $h_{\tilde{f}_t}$ coincide all of them with the identity near $0$ and $1$ and such that, for all $t<1$, 
$\Vert h_{\tilde{f}_t}fh_{\tilde{f}_t}^{-1}-f_t\Vert_{1}<\varepsilon_t$. 
Moreover, Corollary \ref{c2. densité C1} ensures the $\mc{C}^{1}$-continuity of the path $(h_{\tilde{f}_t})_{t\in[0;1[}$.
\hfill$\Box$\medskip

\subsection{Mather invariant}

~~We will now reword Propositions \ref{p1.densité C1} and \ref{p1.densité C1 para} by stressing what makes $h_{g}$ and $id$ be different from each other near $1$.

~~For that, we will introduce a new notion, which will be an equivalent of Mather invariant which appeared in $\mc{C}^2$-context.

\begin{defi}\label{d.densité C1}
Let us consider $f\in D_{\alpha,\beta}$ and $g\in D_f$.
\\
We will call \emph{ Mather invariant of $g$ with respect to $f$}, and we will denote by $\mc{M}_f(g)$, the unique homeomorphism of $[0;1]$ commuting with $f$ and coinciding with $h_g$ on a neighbourhood of $1$.
\end{defi}

\underline{Remark :} Mather invariant is well-defined since $h_g$ commutes with $f$ on a neighbourhood of $1$. By choosing a fundamental domain $I$ included in this neighbourhood of $1$ and by pushing $\restriction{h_g}{I}$ by the dynamic of $f$, one obtains the unique homeomorphism $\restriction{\mc{M}_f(g)}{]0;1[}$ of $]0;1]$ which commutes with $f$ and coincides with $h_g$ on a neighbourhood of $1$. By finally setting $\mc{M}_f(g)(0)=(0)$, one obtains a diffeomorphism of $[0;1]$.

\begin{lem}\label{l1.densité C1}
If $f\in D_{\alpha,\beta}$ and $g,g'\in D_f$ are such that $g'\geqslant g$ and $g'>g$ on a closed interval containing a fundamental domain of $g'$, then $M_f(g')>M_f(g)$.
\end{lem}

\underline{Proof of Lemma \ref{l1.densité C1} :}
Let $a,b\in]0;1[$ be two real numbers such that $g,g'\in D_{f,a,b}$.
By definition of $\mc{M}_f$, it is sufficient to show that $h_{g'}>h_g$ on a fundamental domain $I_0=[x_0;f(x_0)]$ of $f$ included in $[b;1]$ ; in other words, it is sufficient to show that 
$h_{g'}h_g^{-1}>id$ on $h_{g}(I_0)$, 
or also to show that : $g'^ng^{-n} >id$ on $h_{g}(I_0)$ if $n$ is great enough so that
$g^{-n}(x)<a$.
Let us consider $x\in h_{g}(I_0)$.
Since $g'\geqslant g$, we know that $g'^ng^{-n}(x)\geqslant x_0$, and so the suborbit $\{g'^i(g^{-n}(x)\}_{0\leqslant i<n}$ goes in the fundamental domain of 
$g'$ on which $g'>g$.
Consequently, $g'^n(g^{-n}(x))>g^n(g^{-n}(x))>x$.

\hfill$\Box$\medskip

\begin{lem}\label{l1.densité C1 bis}
Let us consider $f\in D_{\alpha,\beta}$ and $(a_t)_{t\in[0;1[},(b_t)_{t\in[0;1[}$ two continuous paths of real numbers of $]0;1[$. Let $(f_t)_{t\in[0;1[}$ be a path of $\mc{C}^1$-diffeomorphisms such that, for all $t<1$, $f_t\in D_{f,a_t,b_t}$. 
Then, for all $I$ compact interval of $]0;1[$, the path
$(\restriction{\mc{M}_{f}(f_t)}{I})_{t\in[0;1[}$
is continuous with respect to the $\mc{C}^1$-topology.
\end{lem}

\underline{Proof of Lemma \ref{l1.densité C1 bis} :}
Let $t\in[0;1[$ be a real number and $\mc{V}$ be a neighbourhood of $t$.
From Lemma \ref{l1.conjugante}, we know that $h_{f_{t'}}$ continuously depends on $t'\in\mc{V}$ on a fundamental domain $I_0$ of $f$ included in
$[\sup\limits_{t'\in\mc{V}}b_{t'};1[$.
As $I$ is compact, there exists an integer $n$ such that at each point $x\in I$, 
$\mc{M}_f(f_{t'})$ is obtained by conjugating $\restriction{h_{f_{t'}}}{I_0}$ by $f$ less 
than $n$ times. Thus it suffices that $\restriction{h_{f_{t'}}}{I_0}$ is sufficiently close to $\restriction{h_{f_{t}}}{I_0}$, in other words that $t'$ is sufficiently close to $t$ for $\mc{M}_f(f_{t'})$ to be close to $\mc{M}_f(f_{t})$ on whole $I$.
\hfill$\Box$\medskip

Proposition \ref{p1.densité C1} is equivalent to the following proposition : 

\begin{prop}\label{p2.densité C1}
If $f\in D_{\alpha,\beta}$, $g\in D_{f}$, and $\mc{U}$ is a $\mc{C}^{1}$-neighbourhood of $g$, then there exists $g'\in D_{\alpha,\beta}\cap\mc{U}$ such that
$\mc{M}_{f}(g')=id$.
\end{prop}

Similarly, one can reword Proposition \ref{p1.densité C1 para} in words of Mather invariant by : 

\begin{prop}\label{p4.densité C1 para}
Let us consider $f\in D_{\alpha,\beta}$ and $(a_t)_{t\in[0;1[}$ ; $(b_t)_{t\in[0;1[}$ two continuous paths of real numbers of $]0;1[$ ; $(\varepsilon_t)_{t\in[0;1[}$ a continuous path of strictly positive real numbers.
Let lastly $(f_t)_{t\in[0;1[}$ be a $\mc{C}^1$-continuous path of diffeomorphisms of $D_{\alpha,\beta}$, such that, for all $t\in[0;1[$, $f_t$ belongs to $D_{f,a_t,b_t}$.
\\
Then there exist two continuous paths $(\tilde{b}_t)_{t\in[0;1[}$ and $(\tilde{c}_t)_{t\in[0;1[}$ of real numbers of $]0;1[$, as well as a $\mc{C}^1$-continuous path $(\tilde{f}_t)_{t\in[0;1[}$ of diffeomorphisms of $D_{f,a_t,\tilde{b}_t}$ such that, for all $t\in[0;1[$, we have :
\begin{itemize}
\item $\mc{M}_f(\tilde{f}_t)=id$ ;
\item $\Vert \tilde{f}_t-f_t\Vert_{1}<\varepsilon_t$.
\end{itemize}
\end{prop}

Thus we will try, in following sections, to cancel the Mather invariant of a diffeomorphism (i.e. to make it coincide with the identity) while remaining in an arbitrarily small neighbourhood of it.

\subsection{To make the unitary conjugates having a fixed point}

In order to prove Proposition \ref{p2.densité C1}, given $f$ a diffeomorphism of $D_{\alpha,\beta}$ and $g$ belonging to $D_f$, we will first make $\mc{M}_{f}(g)$ having fixed points by modifying $g$ by small perturbations.
More precisely, Proposition \ref{p3.densité C1} will enable us to prescribe the point which will become a fixed point of $\mc{M}_f(g)$ after having perturbed $g$ :

\begin{prop}\label{p3.densité C1}
Let $f$ be an element of $D_{\alpha,\beta}$ ; $a,b\in]0;1[$ be two real numbers such that $a<b$ ;
$g$ be an element of $D_{f,a,b}$ ; $\mc{U}$ be a $\mc{C}^{1}$-neighbourhood of $g$ in 
$D_{\alpha,\beta}$ and $p$ be a point of $]0;1[$.
\\
Then there exist $\tilde{b}>b$ ; $\tilde{g}\in \mc{U}\cap D_{f,a,\tilde{b}}$, such that $p$ is a fixed point of $\mc{M}_{f}(\tilde{g})$.
\end{prop}

The parameter version follows :

\begin{prop}\label{p1.transfochemin}
 Let $f$ be a diffeomorphism of $[0;1]$ without fixed point in $]0;1[$. Let $(f_t)_{t\in[0;1]}$ be a $\mc{C}^1$-continuous path of diffeomorphisms of $[0;1]$ such that :
\begin{itemize}
\item $f_0=f$,
\item $f_t$ has no fixed point in $]0,1[$ for all $t<1$, 
\item For all $t<1$, there exist $a_t>0$ and $b_t<1$, depending continuously on $t$, such that $f_t$ coincides with $f$ on $[0,a_t]\cup[b_t,1]$ (i.e. $f_t\in D_{f,a_t,b_t}$) ;
\end{itemize}
Then, for all path $(\varepsilon_t)_{t\in[0;1[}$ of strictly positive real numbers, for all continuous path $(p_{t})_{\in[0;1[}$ of real numbers of $]0;1[$, there exists a $\mc{C}^1$-continuous path $(\tilde{f}_t)_{t\in[0,1[}$ of diffeomorphisms of $[0;1]$ such that : 
\begin{itemize}
\item $\tilde{f}_0=f$,
\item for all $t<1$, there exists $\tilde{b}_t<1$, depending continuously on $t$, such that $\tilde{f}_t$ coincides with $f$ on $[0;a_t]\cup[\tilde{b}_t;1]$,
\item for all $t\in[0;1[$, $\Vert \tilde{f}_{t}-f_{t}\Vert_{1}<\varepsilon_t$,
\item for all $t\in[0;1[$, $p_t$ is a fixed point of $\mc{M}_f(\tilde{f}_t)$.
\end{itemize}
\end{prop}

~~The idea of the proof is to perturb each $f_{t}$ in a continuous way with respect to $t$. So, one will obtain a new path which will also be continuous with respect to $t$, and such that the unitary conjugacies will have a fixed point for all $t$.
For that, we will consider two diffeomorphisms $f_{+}$ and $f_{-}$ going respectively "arbitrarily faster" and "arbitrarily slowlier" than $f$ -- the meaning of these expressions will be rigorously stated in following Lemma \ref{l1.transfochemin} -- , and we will make coincide the $f_{t}$'s with these diffeomorphisms for a sufficiently long time for them to can catch up or lose their lead with respect to $f$.
\\
First, we will show that such diffeomorphisms $f_{+}$ et $f_{-}$ exist :
\\

\begin{lem}\label{l1.transfochemin}
Let $f$ be a $\mc{C}^{1}$-diffeomorphism of $[0;1]$ such that $f>id$.
There exist two $\mc{C}^{1}$-diffeomorphisms $f_{+}$ and $f_{-}$ of $[0;1]$ such that :
\begin{itemize}
\item $f_{+}>f$, $id<f_{-}<f$ ;
\item $Df_{-}(1)=Df_{+}(1)=Df(1)$ ;
\item for all $x\in]0;1[$, for all $n_{0}\in \N$, there exists $k\in\N$ such that
\begin{eqnarray}\label{eq4} (f_{+})^{k}(x)\geqslant f^{n_{0}+k}(x) \mbox{ and } 
 (f_{-})^{n_{0}+k}(x)\leqslant f^{k}(x) .
\end{eqnarray}
Consequently, for all $k'\geqslant k$, the inequalities 
$(f_{+})^{k'}(x)\geqslant f^{n_{0}+k'}(x)$ and $(f_{-})^{n_{0}+k'}(x)\leqslant f^{k'}(x)$
are still satisfied.
\item «~The lead of $f_{+}$ with respect to $f$ and the delay of $f_{-}$ with respect to $f$ are decreasing~» ; i.e. : given $n_{0}$, the smallest whole number $k$ satisfying the last condition is increasing with respect to $x$.
\end{itemize}
\end{lem}

\underline{Notation :} If $f\in D_{\alpha,\beta}$, $x\in]0;1[$ and $n_0\in\N$, we will denote by $k(n_0,x)$ the smallest whole number $k$ satisfying the condition (\ref{eq4}) above.
\\

\underline{Proof of Lemma \ref{l1.transfochemin} :}
Proof of this lemma is quite long and thus is expounded in Annex.

\hfill$\Box$\medskip

~~Let us recall now the definition of translation number of a diffeomorphism with respect to another commuting with it :

\begin{lede}\label{ld.1}
If $f,h$ are two increasing $\mc{C}^{1}$-diffeomorphisms of $[0;1]$ such that $f$ has no  fixed point on $]0;1[$, is greater than the identity and such that $f$ and $h$ commute on $[0;1]$, then we can consider, given $x\in[0;1]$ and $n\in\N$, the whole number 
$m(n)$ defined by :
\[ f^{m(n)}(x)\leqslant h^n(x)<f^{m(n)+1}(x). \]
Then the limit $\lim_{n\to\infty}\dfrac{m(n)}{n}$ exists and is independent of $x$.
We call it \emph{translation number of $h$ with respect to $f$}, and we will denote it by $\tau_f(h)$.
\\ Furthermore, if $(h_t)_{t\in[O;1[}$ is a path of $\mc{C}^1$-diffeomorphisms which varies $\mc{C}^1$-continuously on the compact sets of $]0;1[$, then the translation number $\tau_f(h_t)$ depends continuously on $t$.
\end{lede}

You can consult \cite{BF} to have a more general proof in the case of local homeomorphisms or \cite{N} for a similar proof in the context of rotation numbers.

Let us also introduce the following definition :

\begin{defi}
Let $f$, $g$ be increasing $\mc{C}^{1}$-diffeomorphisms of $[0;1]$, with no fixed point on $]0;1[$, such that $g$ coincides with $f$ on a neighbourhood $\mc{V}_{0}$ of $0$ and on a neighbourhood $\mc{V}_{1}$ of $1$.
Then the \emph{delay of $g$ with respect to $f$} is the integer
$r_{f}(g)=[\vert\tau_f(\mc{M}_f(g))\vert]$, that is the integral part of the absolute value of the number of translation of the Mather invariant of $g$ with respect to $f$.
\end{defi}

\begin{lem}\label{l2.transfochemin}
Let $f$ be an increasing $\mc{C}^{1}$-diffeomorphism of $[0;1]$ with no other fixed point than $0$ and $1$ ; $(a_t)_{t\in[0;1[}$, $(b_t)_{t\in[0;1[}$ be two continuous paths of real numbers of $]0;1[$ and $(g_t)_{t\in[0;1[}$ be a $\mc{C}^1$-continuous path of increasing diffeomorphisms $g$ of $D_{f,a_t,b_t}$.
Then $\tau_f(\mc{M}_f(g))$ depends continuously on $t$, and the delay of $g_t$ with respect to $f$ is upper semi-continuous with respect to $t$.
\end{lem}

\underline{Proof of Lemma \ref{l2.transfochemin} :}
From Lemma \ref{l1.densité C1 bis}, we know that the Mather invariant of $g_t$ with respect to $f$ varies continuously on the compact sets of $[0;1]$ ; thus, from Lemma -Definition \ref{ld.1}, the number of translation $\tau_f(\mc{M}_f(g))$ depends continuously on $t$, and then one concludes thanks to the continuity of the absolute value and the upper semi-continuity  of the integer part.
\hfill$\Box$\medskip

\begin{lem}\label{l3.transfochemin}
Let $f$ be an element of $D_{\alpha, \beta}$, $(r_t)_{t\in[0;1[}$ be an upper semi-continuous collection of integers and $(x_t)_{t\in[0;1[}$ be a continuous path of real numbers of $]0;1[$.
Then the collection $(\ell_t)_{t\in[0;1[}=(f_+^{k(r_t,x_t)}(x_t))_{t\in[0;1[}$ is 
locally upper bounded ; in other words :
for all $t<1$, there exists a strictly positive real number $\varepsilon$ such that
$\sup\limits_{s\in[t-\varepsilon,t+\varepsilon]}\ell_s<1$.
\\
One obtains of course the same result for the collection $(f_-^{k(r_t,x_t)}(x_t))_{t\in[0;1[}$.
\end{lem}

\underline{Proof of Lemma \ref{l3.transfochemin} :}
Let us consider $t\in[0;1[$ and $\varepsilon>0$. If $s\in[t-\varepsilon;t+\varepsilon]$ and $\varepsilon$ is small enough, then by semi-continuity of $r_s$ one can state that $r_s$ varies from $r_t-1$ to $r_t+1$. Thus $k(r_s,x_s)$ is bounded by 
$\max\limits_{s\in[t-\varepsilon,t+\varepsilon]}k(r_t+1,x_s)$.
Since $f_+$ has been defined in Lemma \ref{l1.transfochemin} in such a way that $k(n_0,x)$ increases with respect to $x$, one can conclude that $k(r_s,x_s)$ is bounded by 
$k(r_t+1,x_{t+\varepsilon})$.
Finally, the increasing of $f_+$ ensures us that 
$\ell_s\leqslant f_+^{k(r_t+1,x_{t+\varepsilon})}(x_{t+\varepsilon})$.

\hfill$\Box$\medskip

\begin{lem}\label{l4.transfochemin}
Let $(\ell_t)_{t\in[0;1[}$ be a locally upper bounded collection of real numbers of $[0;1[$, that is to say satisfying the following property :
for all $t<1$, there exists a strictly positive real number $\varepsilon_t$ such that
$\sup\limits_{s\in[t-\varepsilon_t,t+\varepsilon_t]}\ell_s<1$.
\\
Then there exists a continuous path $(d_t)_{t\in[0;1[}$ such that, for all $t<1$, one has :
 $1>d_t>\ell_t$.
\end{lem}

\underline{Proof of Lemma \ref{l4.transfochemin} :}
Each compact set $C$ of $[0;1[$ is covered by the union of the neighbourhoods $[t-\varepsilon_t;t+\varepsilon_t]$, where $t$ belongs to $C$, from which one can extract finitely many neighbourhoods of this kind which will still cover $C$.
Thus, since $\restriction{\ell}{[t-\varepsilon_t;t+\varepsilon_t]}$ is bounded by $M_t$, $\restriction{\ell}{C}$ is bounded by the maximum of these constants.
So $\ell_t$ is bounded on each compact set of $[0;1[$.
 \\
Let $(t_n)_{n\in\N}$ be an increasing sequence converging to $1$.
On each interval $[t_n;t_{n+1}]$, $\ell$ is bounded by $M_n$, so there exists a continuous path $(d_t)_{t\in[0;1[}$ such that, on $[t_n;t_{n+1}]$, one has $d_t>M_n$, and the lemma is proved.
\hfill$\Box$\medskip



\underline{Proof of Proposition \ref{p1.transfochemin} :}
Since $f_t$ and $f_\pm$ have same derivative at $1$, from Lemma \ref{l3.densité C2 bis}, 
one can exhibit a continuous path $(c_t)_{t\in[0;1[}$ of real numbers of $]0;1[$ satisfying : 
for all $t<1$, $c_t>b_t$, $f^{-1}(\dfrac{c_t+1}{2})>b_t$ and 
$\Vert \restriction{f_t}{[c_t;1]}-\restriction{f_{\pm}}{[c_t;1]}\Vert_{1}
<\tilde{\varepsilon}_t=\dfrac{\varepsilon_t}{4M_\Phi}$.
From Lemmas \ref{l2.transfochemin}, \ref{l3.transfochemin} and \ref{l4.transfochemin}, one can by this way exhibit a continuous path $(d_t)_{t\in[0;1[}$ of real numbers of $]0;1[$ such that, for all
$t\in[0;1[$ : $d_t>f_+^{k(r_t+1,f_+(\dfrac{c_t+1}{2}))+1}(\dfrac{c_t+1}{2})$.
For all $t<1$, from Lemma \ref{l2.densité C1 bis}, there exist two diffeomorphisms $f_{t,+}$ and $f_{t,-}$ such that :
\begin{itemize}
\item $f_{t,\pm}$ coincide with $f_t$ on $[0;c_t]\cup[\dfrac{d_t+1}{2};1]$ ;
\item $f_{t,\pm}$ coincide with $f_\pm$ on $[\dfrac{c_t+1}{2};d_t]$ ;
\item $\Vert f_{t,\pm}-f_t\Vert_{1}<\varepsilon_t$.
\end{itemize}
Let us denote $k_t=k(r_t+1,f_+(\frac{c_t+1}{2}))$.
\\

\underline{Claim 1 :}
For all $t<1$, one has the following inequalities on $]d_t;1[$ : 
$\mc{M}_f(f_{t,+})>id$ and $\mc{M}_f(f_{t,-})<id$.

\underline{Proof of Claim 1 :}
Let us consider $t<1$, $x>d_t$, and $n$ a sufficiently great integer so that $f^{-n}(x)<a_t$.
Let $n_1$ be the smallest integer such that $f^{-n_1}(x)<\dfrac{c_t+1}{2}$.
By construction of $d_t$, one can thus decompose $n_1$ as $n_1=r_f(f_t)+1+k_t+1+n'$, where $n'\in\N^*$.
One has then the following equality : 
\begin{eqnarray}\label{eq5}
\mc{M}_f(f_{t,+})(x)=f_{t,+}^{n_1}\mc{M}_f(f_t)f^{-n_1}(x).
\end{eqnarray}

Moreover, $\mc{M}_f(f_t)(f^{-n_1}(x))\geqslant f^{[\tau_{f}(\mc{M}_f(f_t))]}(f^{-n_1}(x))$ by definition of $\tau_{f}(\mc{M}_f(f_t))$.
Thus 
\begin{eqnarray}\label{eq6}
f_{t,+}^{n_1}(\mc{M}_f(f_t)f^{-n_1}(x))\geqslant
f_{t,+}^{n_1} f^{[\tau_{f}(\mc{M}_f(f_t))]}(f^{-n_1}(x)).
\end{eqnarray}
Since $f_{t,+}\geqslant f$ on $[b_t;1]$ and $f^{-n_1}(x)>b_t$, one has also : 
$f_{t,+}^{r_f(f_t)+1} f^{[\tau_{f}(\mc{M}_f(f_t))]}(f^{-n_1}(x))
\geqslant f^{-n_1}(x)$.
So 
$f_{t,+}^{n_1} f^{[\tau_{f}(\mc{M}_f(f_t))]}(f^{-n_1}(x))
\geqslant f_{t,+}^{k_t+1+n'}(f^{-n_1}(x))$.
\\
Then : $f_{t,+}^{k_t+1+n'}f^{-n_1}(x)=
f_{t,+}^{n'}f_{t,+}^{k_t}(f_{t,+}f^{-n_1}(x))$ ;
furthermore $f_{t,+}f^{-n_1}(x)\in[\dfrac{c_t+1}{2};f_+(\dfrac{c_t+1}{2})[$, thus $f_{t,+}$ coincides with $f_+$ on the interval $[f_{t,+}f^{-n_1}(x);f_+^{k_t-1}(f_{t,+}f^{-n_1}(x))$, which enables us to write : 
$f_{t,+}^{k_t+1+n'}f^{-n_1}(x)=f_{t,+}^{n'}f_+^{k_t}(f_{t,+}f^{-n_1}(x))$.
Since $f_{t,+}f^{-n_1}(x)\in[\dfrac{c_t+1}{2};f_+(\dfrac{c_t+1}{2})[$, and by increasing
of $k$ with respect to its second variable, one knows then that 
$f_{t,+}^{k_t+1+n'}f^{-n_1}(x)>f_{t,+}^{n'}f^{r_f(f_t)+1+k_t}(f_{t,+}f^{-n_1}(x))$.
By noticing that $f_{t,+}=f_+$ is greater than $f$ on 
$[\dfrac{c_t+1}{2};f_+(\dfrac{c_t+1}{2})[$, one can conclude ; indeed one obtains :
$f_{t,+}^{k_t+1+n'}f^{-n_1}(x)>f_{t,+}^{n'}f^{r_f(f_t)+1+k_t}f^{-n_1+1}(x)
\geqslant f_{t,+}^{n'}f^{-n'}(x)\geqslant x$.
\\
From these calculations it follows that $\mc{M}_f(f_{t,+})>id$ on $]d_t;1[$.
\\
One can show similarly that $\mc{M}_f(f_{t,-})<id$ on $]d_t;1[$.

\hfill$\Box$\medskip

Thus, $\mc{M}_f(f_+)$ (resp. $\mc{M}_f(f_-)$ ) is strictly greater (resp. smaller) than the identity on at least one fundamental domain of $f$, and is defined elsewhere by its commuting relation with $f$ ; it is then strictly greater 
(resp. smaller) to the identity on whole $]0;1[$ ; 
in particular $\mc{M}_f(f_+)(p_t)>p_t$ (resp. $\mc{M}_f(f_-)(p_t)<p_t$).
\\

~~One considers now, for all $t<1$, the path of diffeomorphisms 
$(f_{t,s})_{s\in[0;1]}$ defined by : $f_{t,s}=sf_{t,+}+(1-s)f_{t,-}$.
\\

\underline{Claim 2 :}
For all $t<1$, there exists a unique parameter $s_t\in[0;1]$ such that
$\mc{M}_f(f_{t,s_t})(p_t)=p_t$.
\\
One can notice straight away that the property $\Vert f_{t,\pm}-f_t\Vert_{1}<\varepsilon_t$,
which has been stated above, implies that, for all $s\in[0;1]$ : 
\[\Vert f_{t,s}-f_t\Vert_{1}<\varepsilon_t.\] 
Moreover, if we define $\tilde{b}_t=\dfrac{d_t+1}{2}$, we have : 
\[f_{t,s}\in D_{f,a_t,\tilde{b}_t}.\]

\underline{Proof of Claim 2 :}
By construction of $f_{t,+}$, this diffeomorphism coincides with $f_+$ on at least one  fundamental domain, and also is strictly greater than $f_{t,-}$ on this domain, since $f_{t,-}$ coincides on it with $f_-$. Consequently, if $s$ and $s'$ are two real numbers of $[0;1]$ such that $s<s'$, then $f_{t,s'}$ is strictly greater than $f_{t,s}$ on at least one fundamental domain of $f_{t,+}$, and as a consequence also on a least one  fundamental domain of $f_{t,s'}$.
One can thus use Lemma \ref{l1.densité C1} in order to state that, if $s<s'$, 
then $\mc{M}_f(f_{t,s})<\mc{M}_f(f_{t,s'})$.
On the other hand, given $t\in[0;1[$, one knows that $\mc{M}_f(f_{t,s})$ depends continuously on $s$ on each compact set of $]0;1[$, and that $\mc{M}_f(f_{t,+})(p_t)>p_t$ and
$\mc{M}_f(f_{t,-})(p_t)<p_t$.
The result follows.

\hfill$\Box$\medskip

\underline{Claim 3 :}
The real number $s_t$ depends continuously on $t$.

\underline{Proof of Claim 3 :}
Given a real number $t\in[0;1[$ and a sequence $(t_n)_{n\in\N}$ converging to $t$, there exists a whole number $N$ such that, if $n$ is great enough, one has :
$f^{-N}(p_{t_n})=f_{t_n,s_{t_n}}^{-N}(p_{t_n})$.
Moreover, by continuity of the path $(p_t)_{t}$ and of the diffeomorphism $f^{_N}$, one knows that 
$f^{-N}(p_{t_n})$ converges to $f^{-N}(p_t)$, which is equal to $f^{-N}_{t,s_t}(p_t)$.
So $f_{t_n,s_{t_n}}^{-N}(p_{t_n})$ converges to $f^{-N}_{t,s_t}(p_t)$, and by uniqueness of the parameter $s$ satisfying this property, one can deduce that $s_{t_n}\limite{n\to\infty}{}s_t$, which proves the continuity of $(s_t)_{t\in[0;1[}$.

\hfill$\Box$\medskip

~~Claim 3 implies the continuity of the path
$(\tilde{f}_t)_{t\in[0;1[}=(f_{t,s_t})_{t\in[0;1[}$, so the proof is now complete.

\hfill$\Box$\medskip

\subsection{To create a degenerate fixed point of $h_g$} \label{s.point fixe dégénéré}

\begin{lem}\label{l4.densité C1}
Let us consider $f\in D_{\alpha,\beta}$ ; $a,b\in]0;1[$ such that $a<b$ ;
$g\in D_{f,a,b}$ such that $\mc{M}_{f}(g)$ has a fixed point $p$ in $]b;1[$,
and $\mc{U}$ a $\mc{C}^{1}$-neighbourhood of $g$ in 
$D_{\alpha,\beta}$
\\
Then there exist $\tilde{b}>b$ and $\tilde{g}\in \mc{U}\cap D_{f,a,\tilde{b}}$ such that $p$ is a fixed point of $\mc{M}_f(\tilde{g})$, with derivative equal to $1$.
\end{lem}

~~This time again, one states a parameter version of this lemma :

\begin{lem}\label{l4.densité C1 bis}
Let $f$ be a diffeomorphism of $D_{\alpha,\beta}$, $(f_t)_{t\in[0;1[}$ be a  $\mc{C}^{1}$-continuous path of diffeomorphisms of $D_{\alpha,\beta}$, $(p_t)_{t\in[0;1[}$ a continuous path of real numbers of $]0;1[$ and $(\varepsilon_t)_{t\in[0;1[}$ a continuous path of strictly positive real numbers, such that :
\begin{itemize}
\item $f_0=f$ ;
\item there exist two continuous paths $(a_t)_{t\in[0;1[}$ and $(b_t)_{t\in[0;1[}$ of real numbers of $]0;1[$ such that, for all $t<1$, one has $f_t\in D_{f,a_t,b_t}$ ;
\item for all $t<1$, one has : $p_t>b_t$ and $p_t$ is a fixed point of $\mc{M}_f(f_t)$.
\end{itemize}
Then there exist a $\mc{C}^1$-continuous path $(\tilde{f}_t)_{t\in[0;1[}$ of $\mc{C}^{1}$-diffeomorphisms of $D_{\alpha,\beta}$, and a continuous path $(\tilde{b}_t)_{t\in[0;1[}$ of real numbers of $]0;1[$ such that :
\begin{itemize}
\item $\tilde{f}_0=f$ ;
\item $\tilde{f}_t\in D_{f,a_t,\tilde{b}_t}$ ;
\item $p_t$ is a fixed point of $\mc{M}_f(\tilde{f}_t)$ with derivative equal to $1$ ;
\item for all $t\in[0;1[$, one has $\Vert \tilde{f}_t-f_t\Vert_{1}<\varepsilon_t$.
\end{itemize}
\end{lem}

Here again, since the proof of Lemma \ref{l4.densité C1 bis} is quite simple, one gives directly the proof of the parameter version.
\\

\underline{Proof of Lemma \ref{l4.densité C1 bis} :}
The idea is to perturb locally, along the orbit of $p_t$, the diffeomorphism
$f_t$ by composition on the right by an affine map which will, fundamental domain
by fundamental domain, make the derivative of $h_f(f_t)$ at these fixed points become closer and closer to $1$. 
For that, one chooses a continuous path of strictly positive real numbers, $(\varepsilon'_t)_{t\in[0;1[}$, such that, for all $t<1$, for all 
$\mc{C}^1$-diffeomorphism $\varphi_t$ of $[0;1]$ satisfying 
$\Vert\varphi_t-id\Vert_{1}<\varepsilon'_t$, 
one has $\Vert f_t\circ\varphi-f_t\Vert_{1}<\tilde{\varepsilon_t}$.
For all $t\in[0;1[$, one can consider the smallest whole number $k_t$ such that $Dh_f(f_t)(p_t)(1-\varepsilon'_t)^{k_t}<1$ if $Dh_f(f_t)(p_t)>1$ 
(resp. such that $Dh_f(f_t)(p_t)(1+\varepsilon'_t)^{k_t}>1$ if $Dh_f(f_t)(p_t)<1$).
One considers thus the composition of $f_t$ by the affine map
$H_{1-\varepsilon'_t}^i(x)=(1-\varepsilon'_t)(x-f^i(p_t))+f^i(p_t)$ 
(resp.$H_{1+\varepsilon'_t}^i(x)=(1+\varepsilon'_t)(x-f^i(p_t))+f^i(p_t)$) on a neighbourhood of the fixed points $f^i(p_t)$ where $i=0,\ldots,k_t-2$, and then by the affine map
$\dfrac{1}{(1-\varepsilon'_t)^{k_t-1}Dh-f(f_t)(p_t)}(x-f^{k_t-2}(p_t))+f^{k_t-2}(p_t)$
(resp. $\dfrac{1}{(1+\varepsilon'_t)^{k_t-1}Dh-f(f_t)(p_t)}(x-f^{k_t-2}(p_t))+f^{k_t-2}(p_t)$)
on a neighbourhood of the fixed point $f^{k_t-2}(p_t)$.
One can notice that, as $p_t$ has been chosen to be greater than $b_t$, the diffeomorphism $h_f(\tilde{f}_t)$ at the point $f^i(p_t)$ is obtained by conjugating the diffeomorphism $H_{1\pm\varepsilon'_t}^{i-1}\circ h_f(\tilde{f}_t)$ by $f$
at the point $f^{i-1}(p_t)$. This ensures the preservation of the improvement given to the derivative of $h_f(\tilde{f_t})$ along the orbit of $p_t$, and thus enables us to conclude that, after having worked as explained above, the derivative of $h_f(\tilde{f_t})$ at the fixed point $f^k(p_t)$ is equal to $1$.
\\
The $\mc{C}^{1}$-diffeomorphism $\tilde{f}_t$ is then given by re-gluing these local perturbations to the initial diffeomorphism $f_t$ following Corollary~\ref{c.transfochemin}.
The continuity of this new path $(\tilde{f}_t)_{t\in[0;1[}$ follows from the continuity of the path of fixed points $(p_t)_{t\in[0;1[}$, of the path $(\varepsilon_t)_{t\in[0;1[}$ and from the  continuity of $h_f(f_t)$ with respect to $t$ on each compact set of $]0;1[$.

\hfill$\Box$\medskip

\subsection{To squash $h_g$ in successive fundamental domains : end of the proof of Theorems~\ref{t.densité C1} and \ref{t1.transfochemin} }

In this section, we end the proof of Proposition \ref{p2.densité C1} by the following proposition :

\begin{prop}\label{p4.densité C1} Let $(\alpha_n)_{n\in\N},(\beta_n)_{n\in\N}$ be two sequences of real numbers of $[0;1]$, $(f_{n})_{n\in\N}$ be a sequence of 
diffeomorphisms of $D_{\alpha_n,\beta_n}$ converging to $id$ with respect to the $\mc{C}^1$-topology and 
$h_{0}$ be an increasing $\mc{C}^{1}$-diffeomorphism of $[0;1]$ such that $Dh_0(0)=1=Dh_0(1)$.
Let $\varepsilon>0$ be a real number.
\\
Then there exists a sequence $(\tilde{f}_{n})_{n\in\N}$, where $\tilde{f}_{n}\in 
D_{\alpha_n,\beta_n}\cap\mc{B}_{f_{n}}(\varepsilon)$ for all $n$, such that the sequence of $\mc{C}^{1}$-diffeomorphisms of $[0;1]$ $(h_{n})_{n\in\N}$, defined by $h_0$ and $h_{n}=\tilde{f}_{n-1}h_{n-1}f_{n-1}^{-1}$ if $n\in\N^{*}$, is stationary, equal to $id$ for all $n$ large enough.
\end{prop}
 
Theorem~\ref{p4.densité C1 para} will be proved by a parameter version of this proposition :
 
\begin{prop}\label{p1.findémo}
Let $(f_{t,n})_{(t,n)\in[0;1[\times\N}$ be a collection of diffeomorphisms of $D_{\alpha_n,\beta_n}$ such that : 
\begin{itemize}
\item for all $n$, $(f_{t,n})_{t\in[0;1[}$ is a $\mc{C}^{1}$-continuous path, 
\item for all $t\in[0;1[$, 
$(f_{t,n})_{n\in\N}$ converges to the identity when $n$ tends to infinity in $\mc{C}^1$-topology.
 \end{itemize}
 Let also $(h_{t,0})_{t\in[0;1[}$ be a continuous path of increasing $\mc{C}^{1}$-diffeomorphisms of $[0;1]$ such that, for all $t\in[0;1[$, one has $Dh_{t,0}(0)=1=Dh_{t,0}(1)$, and lastly let $(\varepsilon_t)_{t\in[0;1[}$ be continuous path of strictly positive real numbers.
 \\
 Then there exists a collection $(\tilde{f}_{t,n})_{(t,n)\in[0;1[\times\N}$, such that $\tilde{f}_{t,n}\in D_{\alpha_n,\beta_n}\cap\mc{B}_{f_{t,n}}(\varepsilon_t)$ for all $(t,n)\in[0;1[\times\N$, and such that the collection of $\mc{C}^{1}$-diffeomorphisms of $[0;1]$ $(h_{t,n})_{n\in\N}$, defined by $h_{t,0}$ and $h_{t,n}=\tilde{f}_{t,n-1}h_{t,n-1}f_{t,n-1}^{-1}$ if $n\in\N^{*}$, is stationary for all $t$, equal to $id$ for all $n$ greater than a whole number $N$ great enough.
 \end{prop}

~~Let us make a technical remark before beginning the proof of this proposition.
 
 \underline{Remark :}
 The graph of the map $\begin{array}{llll} F : & ]0;+\infty] & \longrightarrow & ]0;+\infty] \\
                   & x           & \longmapsto & \left\vert\dfrac{x}{1-x}\right\vert

 \end{array}$
 is given by the below figure.
 In particular, $F$ is increasing on $]0;1[$ and decreasing on $]1;+\infty[$ ; its limits 
are $+\infty$ at $1^{-}$ and at $1^{+}$ ; it has value $0$ at $0$ and tends to $1$ at $+\infty$ ;
its values at $\dfrac{1}{2}$ is $1$ and at $\dfrac{1}{3}$ is $\dfrac{1}{2}$.
\\
\\

\includegraphics[width=300px,height=300px, keepaspectratio=true]{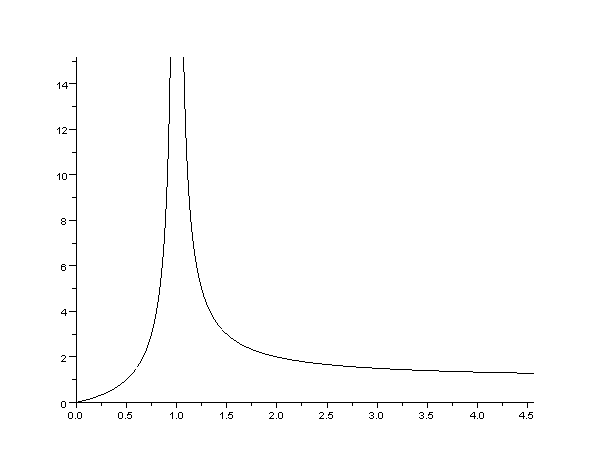}

\underline{Notations :} If $\varepsilon$ is a strictly positive real number, we will denote by $0<$\emph{$a_{\varepsilon}$}$<1<$\emph{$b_{\varepsilon}$}
the two real numbers such that : $F(x)>\dfrac{1}{\varepsilon}$ if $a_{\varepsilon}<x<b_{\varepsilon}$
and $F(x)\leqslant\dfrac{1}{\varepsilon}$ if $x\leqslant a_{\varepsilon}$
or $x\geqslant b_{\varepsilon}$.
\\
\\
 
 In order to expound clearly the reasoning, and as the adaptation to the parameter version does not present any specific difficulty, we show here only how proving Theorem \ref{t.densité C1} by Proposition \ref{p4.densité C1}. 
\\

\underline{Proof of Theorem \ref{t.densité C1} by Proposition \ref{p4.densité C1} :}
 Let us denote by $a,b$ the two real numbers such that $g\in D_{f,a,b}$. Let $p$ be an element of $]0;1[$. From Proposition \ref{p3.densité C1}, there exist $\tilde{b}>b$ and $\tilde{g}\in\mc{U}\cap D_{f,a,\tilde{b}}$ such that $p$ is a fixed point of $\mc{M}_{f}(\tilde{g})$. From Lemma \ref{l4.densité C1}, one can suppose that $D\mc{M}_{f}(\tilde{g})(p)=1$. Let $m$ be a sufficiently great integer so that $f^m(p)>\tilde{b}$.
 One considers the sequence $(f_n)_{n\in\N}$, where $f_n$ is the normalization of 
the diffeomorphism $\restriction{f}{[f^{m+n}(p);f^{m+n+1}(p)]}$ from $[f^{m+n}(p);f^{m+n+1}(p)]$ on its image. One defines also the $\mc{C}^1$-deffeomorphism $h_0$ as being the normalization on $[0;1]$ of $\restriction{\mc{M}_f(g)}{[f^m(p);f^{m+1}(p)]}$. Then this sequence converges uniformly to the identity with respect to the $\mc{C}^1$-topology, and one can apply Proposition \ref{p4.densité C1}. By perturbing $g$ on the successive fundamental domains of $f$ in such a way that $\tilde{f}_n$ is the normalization on $[0;1]$ of the diffeomorphism $\restriction{g}{[f^{m+n}(p);f^{m+n+1}(p)]}$, one will obtain a perturbation$\tilde{g}$ of $g$ on $[0;1]$ such that $\mc{M}_f(\tilde{g})=id$.
\hfill$	\Box$\medskip

\underline{Proof of Proposition \ref{p1.findémo} :}
We will write the diffeomorphisms $\tilde{f}_{t,n}$ in the form of perturbations $f_{t,n}\circ\psi_{t,n}$ of the elements of the sequence $(f_{t,n})_{n}$. 
 By noticing that $h_{t,n+1}$ is then given by $f_{t,n}\circ\psi_{t,n}h_{t,n}f_{t,n}^{-1}$, we can then interpret the transition from $h_{t,n}$ to $h_{t,n+1}$, where $t$ is given, as being a transformation of $h_{t,n}$ into $\psi_{t,n}\circ h_{t,n}$, followed by a conjugacy by $f_{t,n}$.
Considering this, since our aim is to obtain a diffeomorphism $h_{t,n}$ equal to the identity, it appears that the best perturbation $\psi_{t,n}$ would be $h_{t,n}^{-1}$. However, the permitted perturbations have a size which is bounded by $\tilde{\varepsilon}_t$, and this strictly positive real number is defined as being small enough for the condition $\tilde{f}_{t,n}=f_{t,n}\circ\psi_{t,n}\in \mc{B}_{f_{t,n}}(\varepsilon_t)$ to be satisfied. This real number does not depend on $f_{t,n}$, but only on $\varepsilon_t$, and one can choose it in such a way that the path $(\tilde{\varepsilon}_t)_{t\in[0;1[}$ is continuous.

~~Given a whole number $n\in\N$ and a real number $t\in[0;1[$, one defines thus $\psi_{t,n}(x)=x+K_{t,n}\Psi_{t,n}(x)$, 
where $\Psi_{t,n}(x)$ is defined by $h_{t,n}^{-1}(x)=x+\Psi_{t,n}(x)$, and
$K_{t,n}=\inf\left(1;\dfrac{\tilde{\varepsilon}_t}{\max\limits_{x\in [0;1]}\vert D\Psi_{t,n}(x)\vert}\right)=\inf\left(1;\tilde{\varepsilon}_t\cdot\min\limits_{x\in [0;1]}\dfrac{1}{\vert D\Psi_{t,n}(x)\vert}\right)$.
\\
 One can then check that $\Vert \psi_{t,n}-id\Vert_{1}<\tilde{\varepsilon}_t$, and that 
$K_{t,n}$ can be written in the following form : $K_{t,n}=\inf\left(1;\tilde{\varepsilon}_t\cdot\min\limits_{[0;1]}\left\vert\dfrac{Dh_{t,n}}{1-Dh_{t,n}}\right\vert\right)$.
\\
As it was announced, one defines also, if $(t,n)\in[0;1[\times\N$, the diffeomorphism $\tilde{f}_{t,n}$ as being $\tilde{f}_{t,n}=f_{t,n}\circ\psi_{t,n}$, and lastly $h_{t,n+1}$ as being 
$h_{t,n+1}=\tilde{f}_{t,n}h_{t,n}f_{t,n}^{-1}$.
\\
Let us introduce the following operator $P$ :
 \[ \begin{array}{llll}
 P: & \lbrace f_{t,n} ; (t,n)\in[0;1[\times\N\rbrace\times\mc{D}iff^1_+([0;1]) &                           \longrightarrow & \lbrace f_{t,n} ;(t,n)\in[0;1[\times\N\rbrace\times\mc{D}iff^1_+([0;1])\\
    & (f_{t,n},h)                & \longmapsto     & (f_{t,n+1},\tilde{f}_{t,n}hf_{t,n}^{-1})
 \end{array} \] 
 
 In particular, if $(t,n)\in[0;1[\times\N$, one has $P({f}_{t,n},h_{t,n})=({f}_{t,n+1},h_{t,n+1})$.
\\

\underline{Remark :} The operator $P$ is continuous with respect to $t$.
\\
Proof : It follows from the continuity of the paths $({f}_{t,n})_{t\in[0;1[}$ and $(h_{t,n})_{t\in[0;1[}$ where $n$ is given ; the continuity of $(h_{t,n})_{t\in[0;1[}$ being proved by induction. Indeed, 
 $(h_{t,0})_{t\in[0;1[}$ is assumed to be continuous, and if $(h_{t,n})_{t\in[0;1[}$ is continuous, then it goes same for $\Psi_{t,n}$ and for $K_{t,n}$ (from the continuity of the path $(\tilde{\varepsilon}_t)_{t\in[0;1[}$), and consequently for $\psi_{t,n}$. Let us stress now that the continuity of the path $(\tilde{f}_{t,n})_{t\in[0;1[}$ is then proved, and this implies the continuity of $(h_{t,n})_{t}$.
 \hfill$\Box$\medskip
 
 Let us notice now that, if there exists $(t,n)\in[0;1[\times\N$ such that $\min\limits_{[0;1]}\left\vert\frac{Dh_{t,n}}{1-Dh_{t,n}}\right\vert\geqslant\dfrac{1}{\tilde{\varepsilon}_t}$, then $K_{t,n}=1$ and $h_{t,n+1}=id$. Moreover, one can easily check that, if $n$ is such that $h_{t,n}=id$, then $h_{t,k}=id$ for all $k\geqslant n$.
 \\
 So, if we prove the following, the proof will be complete :
\\

\begin{center}
Claim 1 : For all $t\in[0;1[$, there exists $M_t\in\N$ such that $\min\limits_{[0;1]}\left\vert\frac{Dh_{t,M_t}}{1-Dh_{t,M_t}}\right\vert\geqslant \dfrac{1}{\tilde{\varepsilon}_t}$.
\end{center}

The aim of the rest of the proof is hence to prove this result.
 For that, we will first use the convergence of $(f_{t,n})_{n\in\N}$ to the identity in order to clarify the calculations :
 Let $(t,n)$ be an element of $[0;1[\times\N$. Let $\eta_t>0$ be such that $(1+\tilde{\varepsilon}_t)(1-\eta_t)>1+\dfrac{\tilde{\varepsilon}_t}{2}$ and
$(1-\tilde{\varepsilon}_t)(1+\eta_t)<1-\dfrac{\tilde{\varepsilon}_t}{2}$, such that also
$\eta_t<1-\dfrac{2}{3(1+\frac{\tilde{\varepsilon}_t}{2})}$.
There exists $N_t>0$ such that, for all $n\geqslant N_t$, and all $x,y\in[0;1]$, one has 
$1-\eta_t<\dfrac{Df_{t,n}(x)}{Df_{t,n}(y)}<1+\eta_t$.
\\

Let us consider such a whole number $N_t$, and define now $\eta_{N_t}=\eta_t$ ; the  uniform convergence of $f_{t,n}$ to $id$ when $n$ tends to $\infty$
implies also the existence of a sequence $(\eta_{t,n})_{n\geqslant N_t}$ converging to $0$ and such that, for all $n\geqslant N_t$, for all $x,y\in [0;1]$, one has 
$1-\eta_{t,n}<\dfrac{Df_{t,n}(x)}{Df_{t,n}(y)}<1+\eta_{t,n}$.
One can also demand that $\eta_{t,n}\leqslant\eta_{t}$ for all $t<1$ and all $n\geqslant N_t$.
\\We have the following result :
\\

\underline{Claim 2 :}
 For all $t\in[0;1[$, there exists $n\geqslant N_t$ such that $\min\limits_{[0;1]}\left\vert\frac{Dh_{t,n}}{ 1-Dh_{t,n}}\right\vert>1$.
 \\
 \underline{Proof of Claim 2 :}
 Let us consider $x\in[0;1]$ and $n\geqslant N_t$. Assume that $\min\limits_{[0;1]}\left\vert\frac{Dh_{t,n}}{ 1-Dh_{t,n}}\right\vert\leqslant 1$. We can calculate :
 \begin{eqnarray}\label{eq0}
D(\psi_{t,n}\circ h_{t,n})(x) = Dh_{t,n}(x)+K_{t,n}(1-Dh_{t,n}(x)) .
\end{eqnarray}

~~Let us then notice that : 

\begin{enumerate}
\item If $x$ is such that $Dh_{t,n}(x)=1$, then $1-Dh_{t,n}(x)=0$ thus 
the derivative at this point does not change when
$h_{t,n}$ gets composed by $\psi_{t,n}$.
\item If $x$ is such that $Dh_{t,n}(x)<1$ (resp. $Dh_{t,n}(x)>1$), then 
$1-Dh_{t,n}(x)>0$ (resp. $1-Dh_{t,n}(x)<0$), thus 
the derivative becomes greater (resp. smaller) after having been composed by $\psi_{t,n}$.

\item  If $Dh_{t,n}(x)<1$, then 
$D(\psi_{t,n} \circ h_{t,n})(x)$ is also $<1$, and conversely.
\item if $x,y\in[0;1]$ are such that $Dh_{t,n}(x)<Dh_{t,n}(y)$, then 
$D(\psi_{t,n}\circ h_{t,n})(x)<D(\psi_{t,n}\circ h_{t,n})(y)$. 
\end{enumerate}

If $\min\limits_{[0;1]}\left\vert\frac{Dh_{t,n}}{ 1-Dh_{t,n}}\right\vert>1$, then the claim is 
true. We assume then that $\min\limits_{[0;1]}\left\vert\frac{Dh_{t,n}}{ 1-
Dh_{t,n}}\right\vert\leqslant1$. The graph of $F$ indicates that this minimum is reached at $\min\limits_{[0;1]}Dh_{t,n}$ ; in other words : 
$\min\limits_{[0;1]}\left\vert\frac{Dh_{t,n}}{ 1-Dh_{t,n}}\right\vert
=\frac{\min\limits_{[0;1]} Dh_{t,n}}{\vert 1-\min\limits_{[0;1]} Dh_{t,n}\vert}$. On the other hand, from item 4. it follows that
$\mathrm{argmin}D(\psi_{t,n}\circ h_{t,n})=\mathrm{argmin}Dh_{t,n}$. By using these two elements, one can now calculate $\min D(\psi_{t,n}\circ h_{t,n})=(1+\tilde{\varepsilon}_t)\min Dh_{t,n}$,
i.e. : the minimum of $D(\psi_{t,n}\circ h_{t,n})$ and $Dh_{t,n}$ are reached at the same point.
It follows immediately that, for all $x\in[0;1]$, 
\begin{equation}\label{eq1}
D(\psi_{t,n}\circ h_{t,n})(x)\geqslant(1+\tilde{\varepsilon}_t)\min\limits_{[0;1]}Dh_{t,n} .
\end{equation}
So as to know the derivative of $h_{n+1}$, one has still to conjugate $\psi_{t,n}\circ h_{t,n}$ by $f_{t,n}$, which leads to the following :
\begin{eqnarray*}
Dh_{t,n+1}(f_{t,n}(x))= \dfrac{Df_{t,n}(\psi_{t,n}h_{t,n}(x))}{Df_{t,n}(x)}
            \cdot D(\psi_{t,n}h_{t,n})(x) \mbox{ .}
\end{eqnarray*}
From this expression and from \ref{eq1}, and according our choice for $\eta_{t}$, we deduce :
for all $x\in[0;1]$, 
\begin{equation}\label{eq2}
Dh_{t,n+1}(x)>(1+\tilde{\varepsilon}_t)(1-\eta_{t,n})\min\limits_{[0;1]} Dh_{t,n}>
(1+\dfrac{\tilde{\varepsilon}_t}{2})\min\limits_{[0;1]} Dh_{t,n}.
\end{equation}
Let now $k$ be an integer such that 
$(1+\dfrac{\tilde{\varepsilon}_t}{2})^k\min\limits_{[0;1]} Dh_{t,n}>1$.
Two cases can then occur : 
\begin{itemize}
\item Either there exists an integer $n< k'\leqslant k$ such that 
$\min\limits_{[0;1]}\left\vert\frac{Dh_{t,n}}{ 1-Dh_{t,n}}\right\vert>1$ ; in this case our statement is proved ;
\item Or each integer $k'$ being between $n$ and $k$ is suh that $\min\limits_{[0;1]}\left\vert\frac{Dh_{t,n}}{ 1-Dh_{t,n}}\right\vert\leqslant 1$ ; and in this case, from the previous calculation, one obtains : $Dh_{n+k}>(1+\dfrac{\tilde{\varepsilon}_t}{2})^k\min\limits_{[0;1]} Dh_{t,n}>1$,
\end{itemize}
and from the graph of the map $F$ we can conclude that Claim 2 is proved.
\hfill$\Box$\medskip

\underline{Claim 3 :}
Let $t\in[0;1[$ be a real number. If $n\geqslant N_t$ is such that $\dfrac{1}{\tilde{\varepsilon}_t}>\min\left\vert\frac{Dh_{t,n}}{ 1-Dh_{t,n}}\right\vert> 1$, then for all $k\geqslant n$ one has : $\min\left\vert\frac{Dh_{t,k}}{ 1-Dh_{t,k}}\right\vert >1$.

\underline{Proof of Claim 3 :}
We prove this claim by induction.
If $t,n$ are such that $n\geqslant N_t$ and such that $\min\limits_{[0;1]}\left\vert\frac{Dh_{t,n}}{ 1-Dh_{t,n}}\right\vert> 1$, two cases can occur : 
\begin{itemize}
\item if this minimum is reached at a point $x\in[0;1]$ such that $Dh_{t,n}(x)<1$, then, since 
$\dfrac{1}{\tilde{\varepsilon}_t}>\min\limits_{[0;1]}\left\vert\frac{Dh_{t,n}}{ 1-Dh_{t,n}}\right\vert$, 
one can make similar calculations as those made above which lead to the inequality \ref{eq2}. In particular, it follows : 
$\min\limits_{[0;1]} Dh_{t,n+1} > \min\limits_{[0;1]} Dh_{t,n}$, and we deduce, with the help of the graph of $F$, that 
$\min\limits_{[0;1]}\left\vert\frac{Dh_{t,n+1}}{ 1-Dh_{t,n+1}}\right\vert
>\min\limits_{[0;1]}\left\vert\frac{Dh_{t,n}}{ 1-Dh_{t,n}}\right\vert>1$.

\item otherwise, this minimum is reached at a point $x$ suchthat $Dh_{t,n}(x)>1$ ; thus, from the graph of $F$, it is reached at $\max\limits_{[0;1]} Dh_{t,n}$, that is : $\min\limits_{[0;1]}\left\vert\frac{Dh_{t,n+1}}{ 1-Dh_{t,n+1}}\right\vert=\dfrac{\max\limits_{[0;1]} Dh_{t,n}}{\max\limits_{[0;1]} Dh_{t,n}-1}$. 
Furthermore, as previously we can use item 4. to conclude that $\mathrm{argmax}(D\psi_{t,n}\circ h_{t,n})=\mathrm{argmax}Dh_{t,n}$ ; which gives : $\max\limits_{[0;1]} D(\psi_{t,n}\circ h_{t,n})=(1+\tilde{\varepsilon}_t)\max\limits_{[0;1]} Dh_{t,n}$.
Consequently, for all $x\in[0;1]$, one has : 
\begin{eqnarray}\label{eq3}
Dh_{t,n+1}(x)<(1-\tilde{\varepsilon}_t)(1+\eta_{t,n})\max\limits_{[0;1]} Dh_{t,n}
<\max\limits_{[0;1]} Dh_{t,n}.
\end{eqnarray}
\\
On the other hand, the hypothesis 
"$\min\limits_{[0;1]}\vert\frac{Dh_{n}}{1-Dh_{n}}\vert>1$" implies that each $Dh_{t,n}(x)$ is 
strictly greater than $\dfrac{1}{2}$. One can thus write, for all $x\in[0;1]$ :
$D(\psi_{t,n}\circ h_{t,n})(x)=Dh_{t,n}(1-K_{t,n})+K_{t,n}>\dfrac{1}{2}(1-K_{t,n})+K_{t,n}$, and then, using the hypothesis "$\min\limits_{[0;1]}\vert\frac{Dh_{n}}{1-Dh_{n}}\vert>1$", one obtains : $D(\psi_{t,n}\circ h_{t,n})(x)>\dfrac{1}{2}(1+\tilde{\varepsilon}_t)$.
After having conjugated $\psi_{t,n}\circ h_{t,n}$ by $f_{t,n}$, one has :
 $Dh_{t,n+1}(x)>(1-\eta_{t,n})D(\psi_{t,n}\circ h_{t,n})(x)>\dfrac{1}{2}(1+\frac{\tilde{\varepsilon}_t}{2})>\dfrac{1}{2}$.
\\
Each $Dh_{t,n+1}(x)$ is then bounded by $\dfrac{1}{2}$ on one side and by $\max\limits_{[0;1]} Dh_{t,n}>1$ on the other side, and thus $\min\limits_{[0;1]}\vert\frac{Dh_{t,n+1}}{1-Dh_{t,n+1}}\vert>1$.
\end{itemize}

\underline{Claim 4 :} For all $t\in[0;1[$, one has : $Dh_{t,n}\cv 1$.

\underline{Proof of Claim 4 :} From reasoning above, one can assume that $\min\limits_{[0;1]}\vert\frac{Dh_{n}}{1-Dh_{n}}\vert>1$ for all $n$ greater than a whole number $n_{0}$.
Let $t\in[0;1[$ be a real number and $n\geqslant n_0$ be a whole number.
From equality \ref{eq0}, we can then deduce following inequalities :
\begin{itemize}
\item If $Dh_{t,n}(x)<1$, then 
$D(\psi_{t,n}h_{t,n})(x)>Dh_{t,n}(x)+\tilde{\varepsilon}_t(1-Dh_{t,n})$, 
and then $Dh_{t,n+1}(f(x))>(1-\eta_{t,n})(Dh_{t,n}(x)+\tilde{\varepsilon}_t(1-Dh_{t,n}))$ ;
\item If $Dh_{t,n}(x)>1$, then 
$D(\psi_{t,n}h_{t,n})(x)<Dh_{t,n}(x)+\tilde{\varepsilon}_t(1-Dh_{t,n})$, 
and then $Dh_{t,n+1}(f(x))<(1+\eta_{t,n})(Dh_{t,n}(x)+\tilde{\varepsilon}_t(1-Dh_{t,n}))$.
\end{itemize}

Denoting by $d_{t,n}$ the maximal distance from $Dh_{t,n}$ to $1$ and using item 4., one deduces :
\\
$d_{n+1}<\max\limits_{[0;1]}\left(\eta_{t,n}+(1-\eta_{t,n})d_{n}(1-\tilde{\varepsilon}_t)) ; 
\eta_{t,n}+(1+\eta_{t,n})d_{t,n}(1-\tilde{\varepsilon}_t)\right)
=\eta_{t,n}+(1+\eta_{t,n})d_{t,n}(1-\tilde{\varepsilon}_t)<d_{t,n}(1-\dfrac{\tilde{\varepsilon}}{2})+\eta_{t,n}$.
\\
Let us then notice that $\dfrac{\eta_{t,n}}{\frac{\tilde{\varepsilon}_t}{2}}$ is an attracting fixed point of the affine map $x\longmapsto x(1-\dfrac{\tilde{\varepsilon}_t}{2})+\eta_{t,n}$, and that $\dfrac{\eta_{t,n}}{\frac{\tilde{\varepsilon}_t}{2}}$ converges to $0$ when $n$ tends to infinity when $t\in[0;1[$ is given.
Consequently, it goes same for maximal distance $d_{t,n}$ from $Dh_{t,n}$ to $1$ when $n$ tends to infinity, and thus Claim 4 is proved.
\\
One can check easily that, for all $(t,n)\in[0;1[\times\N$, 
one has $D\psi_{t,n}(0)=1=D\psi_{t,n}(1)$, and that it follows from the following inequalities : $D\tilde{f}_{t,n}(0)=Df_{t,n}(0)$
and $D\tilde{f}_{t,n}(1)=Df_{t,n}(1)$.
\hfill$\Box$\medskip

\section{Isotopy by conjugacy to the identity} \label{s.isotopie identité}

~~In this section, we prove Theorem~\ref{t.2}.

~~It is clear that, if for some fixed point $x$ of $f$ one has $Df(x)\neq 1$, then for all $t\in[0;1[$, one has $D(h_{t}fh_{t}^{-1})(x)=Df(x)\neq 1$, and consequently $f$ can not converge to $id$ with respect to the $\mc{C}^{1}$-topology.
\\
We will thus now show that the condition is a sufficient condition.
\\

~~For that, we will first consider the case where $f$ has no other fixed point than $0$ and $1$.

\begin{prop}\label{p1.isotopie à id}
Let $f$ be an increasing diffeomorphism of $[0;1]$ such that $Df(0)=1=Df(1)$, and such that 
$\mathrm{Fix}(f)=\{0;1\}$.
Let $g$ be an increasing $\mc{C}^1$-diffeomorphism of $[0;1]$ without hyperbolic fixed point and such that $(f-id)(g-id)\geqslant 0$ on whole $[0;1]$.
Then there exists an isotopy by conjugacy $(f_{t})_{t\in[0;1]}$ from $f$ to $id$ 
such that, for all $t\in[0;1[$, one has the following upper boundary : $\Vert f_{t}-id\Vert_{1}<2\cdot\Vert f-id\Vert_{1}$.
\end{prop}

\underline{Proof of Proposition \ref{p1.isotopie à id} :}
Let $(\tilde{f}_t)_{t\in[0;1[}$ be the continuous path of increasing $\mc{C}^1$-dif\-feo\-mor\-phisms of $[0;1]$ defined by : $f_t=(1-t)f+tg$. 
One has then  the following upper boundary :
$\Vert f_t-g\Vert_{1}\leqslant\Vert f-g\Vert_{1}$.
Furthermore, if $t\in[0;1[$, then $f_t$ has no fixed point ; $f_0=f$ and $f_t\limite{t\to 1}{}g$.
One can thus apply Theorem \ref{t.transfochemin} for a path $(\varepsilon_t)_{t\in[0;1[}$ of strictly positive real numbers bounded by $\Vert f-g\Vert_{1}$, converging to $0$ when 
$t$ tends to $1$, and by this way obtain the existence of an isotopy by conjugacy from $f$ to $g$, denoted by $(f_t)_{t\in[0;1[}$, such that $\Vert \tilde{f}_t-f_t\Vert_{1}<\varepsilon_t$ if $t<1$.
By simple triangle inequality, one has the wanted control for 
$\Vert f_t-g\Vert_{1}$.
 
\hfill$\Box$\medskip

\underline{Proof of Theorem \ref{t.2} :}
Let us first notice that the set of connex components of 
$[0;1]\backslash \mathrm{Fix}(f)$ is countable, thus we can fix a numbering $C_{1},\ldots,C_{n},\ldots$ of the closures of these components.
\\
For all $n\geqslant 1$, denoting by $\Phi_{n}$ the affine map from 
$C_{n}$ into $[0;1]$, we define $f_{n}=\Phi_{n}\restriction{f}{C_{n}}\Phi_{n}$.
\\
Then : $\Vert f_{n}-id\Vert_{1} \limite{n \to \infty}{} 0$.
\\Indeed, by uniform continuity of $D\restriction{f}{C_{n}}$, and since
$\mathrm{diam}(C_{n})\limite{n \to\infty}{} 0$, one obtains that $D\restriction{f}{C_{n}}$ converges to $D\restriction{f}{C_{n}}(\partial C_n)=1$ on 
$C_{n}$, and consequently $\restriction{f}{C_{n}}$ converges to $id$.
\\
One considers then a sequence $(t_{n})_{n\in\N}$ where $t_{0}=0<t_{1}<t_{2}<\ldots<t_{n}<\ldots$ and $t_{n}\limite{n \to \infty}{} 1$.
\\From Proposition \ref{p1.isotopie à id}, for all $n\geqslant 1$, there exists an isotopy by conjugacy $((f_{n})_{s})_{s\in[0;1]}$ from $f_{n}$ to $id$.
One defines then, for all $t\in[0;1]$, the diffeomorphism $f_{t}$ by :
for all $n\geqslant 1$, 
$\restriction{f_{t}}{C_{n}}=\left\{\begin{array}{ll}
                            \restriction{f}{C_{n}} \mbox{ if } t<t_{n-1}\\
                      \Phi_{n}^{-1}(f_{n})_{\frac{t-t_{n-1}}{1-t_{n-1}}}\Phi_{n}
                        \mbox{ if } t_{n-1}\leqslant t\leqslant 1
\end{array}\right.$, and $f_{t}(1)=1$.
\\
The path $(f_{t})_{t\in[0;1[}$ is then continuous : it follows in each component $C_{n}$ from the continuity of isotopies $((f_{n})_{s})_{s\in[0;1]}$, and, as regards continuity on the neighbourhood of fixed points, from the fact that these isotopies coincide with $\restriction{f}{C_{n}}$ at these points.
\\
The continuity at $t=1$ can be proved by the following :
\\
If $\varepsilon$ is a strictly positive real number, there exists $n_{0}\in\N$ such that, for  
all $n\geqslant n_{0}$, $\Vert \restriction{f}{C_{n}}-id\Vert_{1}<\dfrac{\varepsilon}{2}$. Then, for all $t\geqslant t_{n_{0}-1}$,
one has : $f_{t}$ coincides either with $f$, or with $id$, or with $(f_{n})_{s}$ where
$s\in[0;1[$, depending on the component $C_{n}$ in which we work. In each of these cases, 
from Proposition \ref{p1.isotopie à id} one has $\Vert \restriction{f_t}{C_{n}}-id\Vert_{1}<2\cdot\dfrac{\varepsilon}{2}=\varepsilon$.
 
\hfill$\Box$

\section{Generalization : proof of Theorem~\ref{t1.généralisation}}

\subsection{Signature of a diffeomorphism of $[0;1]$}

~~For instance, we saw that two $\mc{C}^1$-diffeomorphisms $f$ and $g$ of $[0;1]$ with no other fixed point than $0$ and $1$, with same derivatives at $0$ and at $1$ and such that 
$(f-id)(g-id)\geqslant 0$ are isotopic by conjugacy one to the other.
We saw also that each $\mc{C}^1$-diffeomorphism of $[0;1]$ with derivative equal to $1$ at each of its fixed point is isotopic by conjugacy to the identity.

We would like now to group together these two results in a more general statement, as well as distinguish the cases in wich there exists an isotopy by conjugacy from $f$ to $g$
from the cases in which such an isotopy does not exist, but in which it is yet possible to obtain a sequence of conjugates of $f$ converging to $g$.

\begin{defi}\label{d.signature}
Let $f$ be a diffeomorphism of $\mc{D}iff_+^1([0;1])$ without hyperbolic fixed point, except possibly $0$ and/or $1$. A countable and well-ordered set
$(\{C_i\}_{i\in I},\prec)$ endowed with a map $\sigma : \{C_i\}_{i\in I}\longrightarrow \{+,-\}$ will be called \emph{signature of $f$} 
if there exists an increasing, one-to-one map $\begin{array}{lcll} \smash[t]
  \Psi : & \{C_i\}_{i\in I} & \longrightarrow & [0;1]\backslash Fix(f)\\
         & C_i    & \longmapsto     & x_i
\end{array}$ such that :
\begin{itemize}
\item For all $i\in I$, $f(x_i)-x_i$ has same sign as $\sigma(C_i)$ ;
\item If $i,j\in I$ are such that $x_i<x_j$, 
then there exists $k\in I$ such that $x_i<x_k\leqslant x_j$ and 
$(f(x_k)-x_k)(f(x_i)-x_i)<0$ ;
\item For all $x\in[0;1]\backslash Fix(f)$, there exists $i\in I$ such that, for all 
$y\in[x;x_i]$ (non oriented interval), $(f(y)-y)(f(x_i)-x_i)\geqslant 0$.
\end{itemize}
\end{defi}

\begin{prop}\label{p1.généralisation}
For all orientation-preserving $\mc{C}^1$-diffeomorphism $f$ of $[0;1]$, without other hyperbolic fixed point than possibly $0$ and $1$, the signature of $f$ exists and is unique up to an orientation-preserving isomorphism.
\end{prop}

\underline{Proof of Proposition \ref{p1.généralisation} :}
\\
$\bullet$ \underline{Existence :} Given a $\mc{C}^1$-diffeomorphism $f$ of $[0;1]$, we will first precise the meaning of the expression \emph{maximal interval on which the sign of $f$ does not change} that has been used in introduction.

For all $x\in[0;1]$ such that $f(x)\neq x$, one considers the set
$\mc{I}_x$ of all intervals $]a;b[\subset[0,1]$ such that
\begin{itemize}
\item $x\in ]a,b[$,
\item for all $y\in ]a,b[, (f(x)-x)(f(y)-y)\geqslant 0$
\item $a$ does not belong to the interior of ${\mathrm{Fix}(f)}$ with respect to the induced topology on $[a;1]$. In other words, $a$ is neither in the interior of $\mathrm{Fix}(f)$ nor the lower extremity of a connex component of the interior of $\mathrm{Fix}(f)$. 
\item $b$ does not belong to the interior of $\mathrm{Fix}(f)$ with respect to the induced topology on $[b;1]$,
\end{itemize}

One can then check that, for all $x\in[0;1]\backslash\mathrm{Fix}(f)$, the interval $I_x=\bigcup\limits_{I\in \mc{I}_x}I$ belongs to the set $\mc{I}_x$ : it is thus the maximal element of this set with respect to the inclusion relation.
\\
One considers now the set $\{I_x ; x\in[0;1]\backslash\mathrm{Fix}(f)\}$. 
Let us stress that, if $x$ and $x'$ are in the same connex component of $[0;1]\backslash\mathrm{Fix}(f)$, then $I_x=I_{x'}$ ; it follows in particular that 
this set is countable.
We will thus denote from now on :
\begin{eqnarray*}
\{I_x ; x\in[0;1]\backslash\mathrm{Fix}(f)\} & = & \{C_i\}_{i\in I}\\
                                             & = & \{I_x\}_{x\in I'},
\end{eqnarray*}
where $I$ is a countable set, and $I'$ is a countable set included in $[0;1]$.
\\
One defines now the map $\Phi$ by :
$\begin{array}{llll}
\Phi : & \{C_i\}_{i\in I} & \longrightarrow & \{+;-\}\\
       &          I_x=C_i & \longmapsto     & \mbox{ sign of }f(x)-x ,
\end{array}$

and the map $\Psi$ by :
$\begin{array}{llll}
\Psi : & \{C_i\}_{i\in I} & \longrightarrow & [0;1]\\
       &          I_x=C_i & \longmapsto     & x. 
\end{array}$

One notices then that, if $x,x'\in I'$ and $x<x'$, then $I_x\cap I_{x'}=\varnothing$.
So, $\{I_x\}_{x\in I'}$ is well-ordered, in the same order than the real numbers $x\in I'$, which enables us to conclude as regards increasing and injectivity of $\Psi$.

~~By construction, the sign of $\Phi(I_x)$ where $x\in I'$ is the one of $f(\Psi(I_x))-\Psi(I_x)$.

Let $x,x'$ be elements of $I'$ such that $x<x'$ and such that $(f(x)-x)(f(x')-x')>0$.
Assume that, for all $y\in[x;x']$, one has $(f(x)-x)(f(y)-y)\geqslant 0$.
Then the set $I_x\cup I_{x'}$ should belong to $\mc{I}_x$, which contradicts the maximality of $I_x$ in $\mc{I}_x$.
Therefore there exists $y\in[x;x']$ such that $(f(y)-y)(f(x)-x)<0$. In particular $y$ is not a fixed point, thus there exists $x''\in I'$ such that $I_y=I_{x''}$.
By increasing and injectivity of $\Psi$, from $x<y<x'$ one deduces $I_x<I_y<I_{x'}$, and then
$x<x''<x'$, and, since $x''\in I_y$, one has also : $(f(x'')-x'')(f(x)-x)<0$.

~~If $x\in[0;1]\backslash\mathrm{Fix}(f)$, then there exists $\tilde{x}\in I'$ such that $x\in I_x=I_{\tilde{x}}$. For all $y\in[x;\tilde{x}]$, one has thus $y\in I_{\tilde{x}}$, so
$(f(y)-y)(f(\tilde{x}-\tilde{x})\geqslant 0$.

$\bullet$ \underline{Uniqueness :}
Let us asssume that $((C=\{C_i\}_{i\in I},\prec),\sigma)$ and $((C'=\{C'_i\}_{i\in I'},\prec'),\sigma')$ are two signatures of an only $\mc{C}^1$-diffeomorphism $f$ of $[0;1]$.
Let us consider $i\in I$ and $x_i\in[0;1]\backslash\mathrm{Fix}(f)$ be the image of $i$
by the map $\Psi$ defined in Definition \ref{d.signature}. Then, since $((C'=\{C'_i\}_{i\in I'},\prec'),\sigma')$ is a signature of $f$, there exists a unique $x'_i$ in the image of $\Psi'$ such that the sign of $f$  does not change on the whole non-oriented interval $[x_i;x'_i]$. We will denote this real number by $\varphi(x_i)$. We define then $\phi(i)=\Psi'^{-1}(x'_i)$. In other words, $\phi$ is the map from $I$ to $I'$ defined by : $\phi=\Psi'^{-1}\varphi\Psi$.
What we will show now is that $\phi$ is isomorphic and order-preserving.
\begin{itemize}
\item $\phi$ is injective : It follows directly from the same properties for the maps $\psi$ and $\Psi'$.
\item $\phi$ is onto : Let us consider $C'_i$ in $C'$. One considers then $x'_i=\Psi'(C'-i)$,
and $x_i$ in the image of $\Psi$ such that the sign of $f$ does not change on 
the non-oriented interval $[x_i;x'_i]$. Then $x'_i$ satisfies the same property as $\varphi(x_i)$, and thus $x'_i=\varphi(x_i)$ ; that enables us to conclude that $C'_i$ is in the image of $\phi$.
\end{itemize}

\hfill$\Box$\medskip

\begin{defi}
From now on, the signature of $f$ \emph{in meaning $A$} will refer to the signature of $f$ as constructed in the proof of the existence of a signature above.
\end{defi}


\subsection{Proof of Theorem~\ref{t1.généralisation} :}

~~First, let us state some properties concerning isotopies by conjugacy :
\\

\underline{Notation :}
If $f,g$ are two orientation-preserving $\mc{C}^1$-diffeomorphisms of $[0;1]$, one will write $f\blacktriangleright g$ if there exists an isotopy by conjugacy from $f$ to $g$.

\begin{prop}\label{p2.généralisation}
The relation $\blacktriangleright$ is reflexive and transitive.
\end{prop}

\underline{Proof of Proposition \ref{p2.généralisation} :}
\\
One immediately checks that $\blacktriangleright$ is reflexive.
\\
Concerning the transitivity, one assumes that $(H_t)_{t\in[0;1[}$ is an isotopy by conjugacy from an orientation-preserving $\mc{C}^1$-diffeomorphism $f_0$ of $[0;1]$
to another, which we will denote by $f$.
One assumes also that $(h_t fh_t^{-1})_{t\in[0;1[}$ is an isotopy by conjugacy from $f$ to a diffeomorphism $g$ of $\mc{D}iff^1_+([0;1])$, and one considers a continuous path of strictly positive real numbers $(\varepsilon_t)_{t\in[0;1[}$.
\\

~~Denoting by $\eta_t>0$ the greatest constant such that, 
if $\Vert \varphi-\psi\Vert_{1}<\eta_t$, 
then $\Vert h_t\varphi h_t^{-1}-h_t\psi h_t^{-1}\Vert_{1}<\varepsilon_t$, the continuity of the path $(h_t)_{t\in[0;1[}$ ensures that the collection 
$(\eta_t)_{t\in[0;1[}$ is locally bounded ; more precisely : 
for all $t\in[0;1[$, 
there exists a neighbourhood $\mc{V}$ of $t$ and $\tilde\eta_{t,\mc{V}}>0$ such that, 
for all $t'\in\mc{V}$, one has : $\tilde\eta_{t,\mc{V}}<\eta_{t'}$.
\\
On the other hand, to each given $t_0\in[0;1[$, by convergence of $H_t$ to $f$, one can associate the smallest real number $\tilde{r}_{t_0}$ of $[0;1[$ such that 
$\Vert f- H_{\tilde{r}_{t_0}}\Vert_{1}<\eta_{t_0}$.
Since the collection $(\eta_t)_{t\in[0;1[}$ is locally bounded, one has a similar property for the collection $(\tilde{r}_t)_{t\in[0;1[}$, that is to say :
for all $t\in[0;1[$, there exists a neighbourhood $\mc{V}$ of $t$ and a real number 
$0\leqslant r_{t,\mc{V}}<1$ such that, 
for all $t'\in\mc{V}$, $\tilde{r}_t<r_{t,\mc{V}}$.
One can therefore use Lemma \ref{l4.transfochemin} which ensures the existence of a continuous path $(r(t))_{t\in[0;1[}$ of real numbers of $[0;1[$ satisfying :
for all $t\in[0;1[$, $\Vert H_{r(t)}-f\Vert_{1}<\eta_t$ ; and consequently also :
$\Vert h_tH_{r(t)}h_t^{-1}-h_tfh_t^{-1}\Vert_{1}<\varepsilon_t$.
By choosing the path $(\varepsilon_t)_{t\in[0;1[}$ in such a way that
$\varepsilon_t\limite{t\to 1}{}0$, and by making tend $t$ to $1$, one obtains the convergence with respect to the $\mc{C}^1$-topology of $h_tH_{r(t)}h_t^{-1}$ to $g$, and consequently one obtains an isotopy by conjugacy from $f_0$ to $g$.
\hfill$\Box$\medskip

\begin{prop}\label{p3.généralisation}
Let $f,g$ be two increasing $\mc{C}^1$-diffeomorphisms of $[0;1]$ without hyperbolic fixed point, and $\mc{B}=\{B_i\}_{i\in I}$ be the set of connex components of 
$[0;1]\backslash\mathrm{Fix}(f)$. One assumes that, for all $i\in I$, one has 
$(\restriction{g}{B_i}-id)(\restriction{f}{B_i}-id)\geqslant 0$, and that $g$ coincides with the identity outside $\bigcup\limits_{i\in I}B_i$.
\\
Then there exists an isotopy by conjugacy from $f$ to $g$.
\end{prop}

\underline{Proof of Proposition \ref{p3.généralisation} :}
This proof is similar to the one of Theorem \ref{t.2}, and can easily be completed by this way by using Proposition \ref{p1.isotopie à id}.
\hfill$\Box$\medskip

\underline{Proof of Theorem \ref{t1.généralisation} (necessary conditions) :}

\begin{enumerate}
\item ~~First, we assume that there exists an isotopy from $f$ to $g$, denoted by $(h_tfh_t^{-1})_{t\in[0;1[}$, where $(h_t)_{t\in[0;1[}$ is continuous. 
One defines $\Psi$ (resp. $\Psi'$) as being the increasing and injective map $\{C_i\}_{i\in I}\longrightarrow [0;1]\backslash Fix(f)$ 
(resp. $\{C'_i\}_{i\in I'}\longrightarrow [0;1]\backslash Fix(g)$) satisfying the conditions listed in the definition of the signature, 
and, if $i\in I$ (resp. $i\in I'$), one defines $x_i=\Psi(C_i)$ (resp. $x'_i=\Psi'(C'_i)$).

Notice then that, if $x'_i\in \Psi'(C')$, then, by convergence of $h_tfh_t^{-1}$ to $g$, there exists $T>0$ such that, for all $t\geqslant T$, one has :
$(h_tfh_t^{-1}(x'_i)-x_i)(g(x_i)-x_i)>0$. So, $f-id$ has at the point $h_t^{-1}(x'_i)$ same sign as $g$ at the point $x'_i$ for all $t\geqslant T$. By definition of the signature of $f$, there exists then an unique $x_{\varphi(i)}\in\Psi(C)$ such that $f-id$ has constant sign on the whole non-oriented interval $[h_T^{-1}(x'_i);x_{\varphi(i)}]$. 
The continuity of $t\longmapsto h_t^{-1}(x'_i)$ ensures that $f-id$ has constant sign on each non-oriented interval of the kind $[h_t^{-1}(x'_i) ; x_{\varphi(i)}]$, where $t\geqslant T$ (and from that we can ensure in particular, that $\varphi$ is well-defined : this map does not depend on the chosen real number $t$, provided that it satisfies the mentioned conditions).
One defines then by $\Phi(C'_i)=\Psi^{-1}(x_{\varphi(i)})$, the map $\Phi$ : $C'\longrightarrow C$, and one shows that :
\begin{itemize}
\item $\Phi$ preserves the signs :
the explanation above shows that $\sigma(\Phi(C'_i))=\sigma'(C'_i)$ for all $i\in I'$.
\item $\Phi$ is non-decreasing : 
Let us consider $C'_i,C'_k\in C'$ such that $C'_i\prec C'_k$ : one has then $x'_i<x'_k$. 
If $g-id$ had same sign at $x'_i$ as at $x'_k$, then, by definition of the signature of $g$, it should exist $x'_l\in \Psi'(C')$, $x'_i<x'_l<x'_k$ such that $(g(x'_i)-x_i)(g(x'_l)-x'_l)<0$. 
Thus it is sufficient to consider the case where $(g(x'_i)-x_i)(g(x'_k)-x'_k)<0$.
One has also, for each $t\in[0;1[$ greater than a sufficiently great real number $T$ of $[0;1]$ : $h_t^{-1}(x'_i)<h_t^{-1}(x'_k)$, 
and $(f(h_t^{-1}(x'_i))-h_t^{-1}(x'_i))(f(h_t^{-1}(x'_k))-h_t^{-1}(x'_k))<0$.
Let us assume now that $x_{\varphi(k)}\leqslant x_{\varphi(i)}$.
Thre cases can then occur :
\\
If $x_{\varphi(k)}\leqslant h_t^{-1}(x'_i)$, then $f-id$ must have constant sign on
the interval $[x_{\varphi(k)};h_t^{-1}(x'_k)]$ ; yet $h_t^{-1}(x'_i)$ belongs to this interval ; that leads us to a contradiction.
\\
If $x_{\varphi(i)}\geqslant h_t^{-1}(x'_k)$, then $f-id$ must have constant sign on the interval $[h_t^{-1}(x_i);x_{\varphi(i)}]$, which contradicts the fact that $h_t^{-1}(x'_k)$ belongs to this interval.
\\
If $h_t^{-1}(x_i)<x_{\varphi(k)}\leqslant x_{\varphi(i)}<h_t^{-1}(x'_k)$, then, since $f-id$ has to have constant signon one hand on $[h_t^{-1}(x'_i);x_{\varphi(i)}]$ and 
on the other hand on $[x_{\varphi(k)};h_t^{-1}(x'_k)]$, one has : 
$f-id$ has constant sign on $[h_t^{-1}(x'_i);h_t^{-1}(x'_k)]$.
This time again, this contradicts the fact that $f-id$ has not the same sign at $h_t^{-1}(x'_i)$ as at $h_t^{-1}(x'_k)$.
\\
Thus $x_{\varphi(i)}<x_{\varphi(k)}$, and then the increasing of $\Psi^{-1}$ enables us to conclude.
\item $\Phi$ is one-to-one : it follows directly from its increasing.
\end{itemize}

\item ~~Let $\{C'_i\}_{i\in\crochetl 1;N\crochetr}$ be a finite subset of $\{C'_i\}_{i\in I'}$. One assumes now that there exists a sequence of conjugates of $f$, denoted by
$(h_nfh_n^{-1})_{n\in\N}$, converging to $g$ when $n$ tends to $\infty$.
Let $m$ be a sufficiently great whole number so that, for all $i\in\crochetl 1;N\crochetr$,
$h_nfh_n^{-1}(x'_i)-x'_i$ has same sign as $g(x'_i)-x'_i$.
Then $f-id$ has same sign at the point $h_m^{-1}(x'_i)$ as $g-id$ has at the point $x'_i$.
By definition of the signature of $f$, 
for all $i\in\crochetl 1;N\crochetr$, there exists $x_{\varphi(i)}\in\Psi(C)$ such that $f-id$ has constant sign on the non-oriented interval $[x_{\varphi(i)};h_m^{-1}(x'_i)]$. 
On defines then $\Phi(C'_i)=\Psi^{-1}(x_{\varphi(i))}$ and
one shows, similarly to case 1., that $\Phi$ is an increasing and sign-preserving one-to-one map.
Here, contrary to the case where we had an isotopy by conjugacy from $f$ to $g$, 
$\varphi$ depends on the chosen whole number $m$, that is why the proof would not work if we had not restricted it to a finite subset of $C'$.
\end{enumerate}

\begin{prop}\label{p4.généralisation}
Let $f$ and $g$ be two increasing $\mc{C}^1$-diffeomorphisms of $[0;1]$ without hyperbolic fixed point.
Let us denote by $((C=\{C_i\}_{i\in I},\prec),\sigma)$ and $((C'=\{C'_i\}_{i\in I'},\prec'),\sigma')$ the respective signatures of $f$ and $g$.
\\
Then :
\begin{enumerate}
\item There exists an isotopy by conjugacy from $f$ to $g$ if and only if there exists a one-to-one and order-preserving map $\Phi : C'\longrightarrow C$ such that, for all $i\in I'$, $\sigma(C'_i)=\sigma'(\Phi(C'_i))$.
\item There exists a sequence of conjugates of $f$ converging to $g$ if and only if, for all finite subset $J'$ of $I'$, there exists a one-to-one and order-preserving map $\Phi : \{C'_i\}_{i\in J'}\longrightarrow C$ such that, for all $i\in J'$, $\sigma(C'_i)=\sigma'(\Phi(C'_i))$.
\end{enumerate}
\end{prop}

\underline{End of the proof of Theorem \ref{t1.généralisation} by Proposition \ref{p4.généralisation} :}

\begin{enumerate}
\item One considers the signature in meaning $A$ of $f$, $((\tilde{C}=\{\tilde{C}_i\}_{i\in I},\tilde{\prec}),\tilde{\sigma})$. 
In each $\tilde{C}_i$, one chooses a point $x_i$ which is not a fixed point of $f$, and one associates to it the connex component $C_i$ of $[0;1]\backslash\mathrm{Fix}(f)$ to which it belongs. One defines $\prec$ by : if $i,j\in I$ , then $C_i\prec C_j$ if and only if $C'_i\prec' C'_j$, and the map $\sigma$ by : if $,j\in I$, then $;\sigma(C_i)=\tilde{\sigma}(\tilde{C}_i)$. Then we choose the set 
$((\{C_i\}_{i\in I},\prec),\sigma)$ as being the description of the signature de $f$, and we will call such a description of the signature of a diffeomorphism of $\mc{D}iff^1_+([0;1])$ : signature \emph{in meaning $B$}.
We will denote by $((C',\prec'),\sigma')$ the signature in meaning $B$ of $g$.
\\
Assume that there exists a one-to-one and increasing map $\Phi$ : $C'\longrightarrow C$ such that, for all $i\in I'$, $\sigma(C'_i)=\sigma'(\Phi(C'_i))$.
If $Df(0)=Dg(0)\neq 1$, then one can consider the smallest element of the ordered set $C'$, which we will denote by $]0;a_g[$, and then $\Phi(]0;a_g[)$ is the smallest element of $C$, which we will denote by $]0;a_f[$.
Otherwise, one has $Df(0)=Dg(0)=1$ and one defines $a_f=a_g=0$.
\\
Similarly, if $Df(1)=Dg(1)\neq 1$ then we can consider the greatest element of $C'$, which we will denote by $]b_g;1[$, and then $\Phi(]b_g;1[)$ is the smallest element of $C$, which we will denote by $]b_f;1[$.
Here again, if on the contrary $Df(1)=Dg(1)=1$, then one defines $b_f=b_g=1$.
\\
Then $\restriction{f}{[a_f;b_f]}$ and $\restriction{g}{[a_g;b_g]}$ are increasing $\mc{C}^1$-diffeomorphisms, without hyperbolic fixed point.
Moreover, $\Phi$ is a one-to-one and increasing map from $C'\backslash (]0;a_g[\cup]b_g;1[)$
to $C\backslash (]0;a_f[\cup]b_f;1[)$, compatible with the signs $\sigma'$ and $\sigma$ on these sets.
By if necessary normalize these two diffeomorphisms $\restriction{f}{[a_f;b_f]}$ and $\restriction{g}{[a_g;b_g]}$ so as to obtain diffeomorphisms defined on the interval $[0;1]$, one can thus use Proposition \ref{p4.généralisation} so as to obtain an isotopy by conjugacy from one to the other.
Let now $h$ be an increasing $\mc{C}^1$-diffeomorphism of $[0;1]$ such that $h(a_f)=a_g$, 
$h(b_f)=b_g$ and which is an affine map on $[a_f;b_f]$. We define $g_1=hgh^{-1}$. Then there exists an 
isotopy by conjugacy from $\restriction{f}{[a_f;b_f]}$ to $\restriction{g_1}{[a_f;b_f]}$.
\\
On $]0;a_f[$, $g_1-id$ has sign $\sigma'([0;a_g])=\sigma([0;a_f])$, and, on $]b_f;1[$,
$g_1-id$ has sign $\sigma'([b_g;1])=\sigma([b_f;1])$. Moreover $g_1-id$ has no fixed point on these intervals ; from Theorem \ref{t.4} there exists an isotopy by conjugacy from $\restriction{f}{[0;a_f]}$ to $\restriction{g_1}{[0;a_f]}$, as well as from 
$\restriction{f}{[b_f;1]}$ to $\restriction{g_1}{[b_f;1]}$.
From these three isotopies by conjugacy, by following the same method as in proof of Theorem \ref{t.2}, one can then construct an isotopy by conjugacy from $f$ to $g_1$, and the transitivity of the relation $\blacktriangleright$ 
enables us to conclude.
\item The proof in the case of a sequence of conjugates follows exactly the same scheme.
\end{enumerate}

\hfill$\Box$\medskip

\begin{prop}\label{p5.généralisation}
Let $f,g$ be two increasing $\mc{C}^1$-diffeomorphisms of $[0;1]$ such that :
\begin{itemize}
\item $(f-id)(g-id)\geqslant 0$, 
\item $Df(0)=Dg(0)=Df(1)=Dg(1)=1$, 
\item $g$ has no fixed point in $]0;1[$,
\item if $x,y\notin\mathrm{Fix}(f)$, then, for all $x\leqslant z\leqslant y$, 
one has $z\notin\mathrm{Fix}(f)$.
\end{itemize}
Then there exists an isotopy by conjugacy from $f$ to $g$, denoted by $(h_tfh_t^{-1})_{t\in[0;1[}$, and a constant $C>0$ such that, for all $t\in[0;1[$, one has : $\Vert h_tfh_t^{-1}-g\Vert_{1}<C\cdot\Vert f-g\Vert_{1}$.
\end{prop}

\underline{Proof of Proposition \ref{p5.généralisation} :}
\\
One proves the result in the case where $f$ has an interval of fixed points $[0;a]$ and has no fixed point in $]a;1[$. The other cases can be solved similarly.
\\
Let $(a_t)_{t\in[0;1[}$ be a continuous path with $a_0=a$ and $a_t\limite{t\to 1}{}0$, and $(h_t)_{t\in[0;1[}$ be a continuous path of diffeomorphisms of $\mc{D}iff^1_{+}([0;1])$, 
such that $h_t(a)=a_t$ and such that $\{\Vert h_t\Vert_{1} ; t\in[0;1[\}$
is bounded (for example one can choose a strictly positive real number $\delta<\min\left(\dfrac{1-a}{2},\dfrac{a}{2}\right)$ and choose $h_t$ coinciding with the affine map from $[0;a]$ into $[0;a_t]$ on $[0;a-\delta]$ and with the affine map from $[a;1]$ into $[a_t;1]$ on 
$[a+\delta;1]$, and then smooth the map on $[a-\delta;a+\delta]$).
\\
Then $h_tfh_t^{-1}$ converges to $f_1$, where $f_1$ has no fixed point on $]0;1[$ and is such that 
$(f-id)(f_1-id)\geqslant 0$ on $[0;1]$.
Then, by using Theorem \ref{t.4} so as to obtain an isotopy by conjugacy from $f_1$ to $g$, and then Proposition \ref{p2.généralisation} so as to put the two isotopies end to end we obtain an isotopy by conjugacy from $f$ to $g$.
\hfill$\Box$\medskip

~~By a similar method to the one presented in the proof of Theorem \ref{t.2}, one deduces from all that the following corollary :

\begin{coro}\label{c.généralisation}
Let $f,g$ be two increasing $\mc{C}^1$-diffeomorphisms of $[0;1]$ without hyperbolic fixed point, and $\mc{B}=\{B_i\}_{i\in I}$ be the set of connex components of 
$[0;1]\backslash\mathrm{Fix}(g)$. One assumes that, for all $i\in I$, 
$(\restriction{g}{B_i}-id)(\restriction{f}{B_i}-id)\geqslant 0$, and that $f$ coincides with the identity outside $\bigcup\limits_{i\in I}B_i$.
\\
Then there exists an isotopy by conjugacy from $f$ to $g$.
\end{coro}

\underline{Proof of Proposition \ref{p4.généralisation} :}
Let us consider the set $((\{C'_i\}_{i\in I'},\prec'),\sigma')$ as being the signature of $g$ in meaning $A$. The set $((\{C_i\}_{i\in I},\prec),\sigma)$, will be considered as being the signature of $f$ in meaning $B$.

Let us consider now a strictly positive real number $\varepsilon$, and let us denote by $J_{\varepsilon}$ the finite subset of $I'$ defined by :
$i\in J_{\varepsilon}$ if there exists $x\in C'_i$ 
such that $\max(\vert g(x)-x\vert,\vert Dg(x)-1\vert)\geqslant\varepsilon$,
so that one has :
$\Vert \restriction{g}{\,\complement\bigcup\limits_{i\in J_n} C'_i}-
\restriction{id}{\,\complement\bigcup\limits_{i\in J_{\varepsilon}} C'_i}\Vert_{1}<\varepsilon$.
\\
Let us denote by $\Phi_{\varepsilon}$ the one-to-one map from $\{C'_i\}_{i\in J_{\varepsilon}}$ to $\{C_i\}_{i\in I}$, 
and let us define the increasing $\mc{C}^1$-diffeomorphism $f_0$ of $[0;1]$
as follows : 

\begin{itemize}
\item on each $\Phi_\varepsilon(C'_i)=]a_i;b_i[$, $f_0$ coincides with the identity on $]a_i;\frac{2}{3}a_i+\frac{1}{3}b_i]\cup[\frac{1}{3}a_i+\frac{2}{3}b_i;b_i[$ ; has no fixed point on $]\frac{2}{3}a_i+\frac{1}{3}b_i;\frac{1}{3}a-i+\frac{2}{3}b_i[$ ; and on this latter interval, $f_0-id$ has same sign as $f-id$ ;
\item $f_0$ coincides with the identity elsewhere.
\end{itemize}
By using again Proposition \ref{p3.généralisation}, one can claim that there exists an isotopy by conjugacy from $f$ to $f_0$.

~~Let us define now the diffeomorphism $f_1$ of $\mc{D}iff^1_+([0;1])$ by :
\begin{itemize}
\item for all $i\in J_{\varepsilon}$, $\restriction{f_1}{C'_i}$ is conjugate to 
$\restriction{f_0}{\Phi_\varepsilon(C'_i)}$ (by an increasing and affine $\mc{C}^1$-diffeomorphism, denoted by $\tilde{h}_i$) ; 
\item $f_1$ coincides with the identity outside $\bigcup\limits_{i\in J_{\varepsilon}}C'_i$ ;
\end{itemize}

Remark : The signature of $f_1$ in meaning $B$ is $(\{C'_i\}_{i\in J_{\varepsilon}}, \sigma')$.

Then : 
\\

\underline{Claim :}
There exists an isotopy by conjugacy from $f_0$ to $f_1$.

Proof of the Claim :
One defines $h_0=id$, and, if $i\geqslant 1$, one defines :
\begin{itemize}
\item $\varphi(i)$ the whole number stricly lower than $i$ satisfying 
$a_{\varphi(i)}=\max\{a_j ; j<i \mbox{ and } a_j< a_i\}$ ; if this set is empty, then one defines $\varphi(i)=0$ ;
\item $\psi(i)$ the whole number stricly lower than $i$ satisfying 
$a_{\psi(i)}=\min\{a_j ; j<i \mbox{ and } a_j> a_i\}$ ; if this set is empty, then one defines $\psi(i)=1$ ;
\item a conjugacy $h_i\in\mc{D}iff^1_+([0;1])$ satisfying :
$\restriction{h_i}{[\frac{2}{3}a'_i+\frac{1}{3}b'_i;\frac{1}{3}a'_i+\frac{2}{3}b'_i]}$ coincides
with $\tilde{h}_i$ ; 
\\
on $[b'_{\varphi(i)};a'_i]$, $h_i$ coincides with the affine map
from $[\frac{1}{3}a'_{\varphi(i)}+\frac{2}{3}b'_{\varphi(i)};\frac{2}{3}a'_i+\frac{1}{3}b'_i]$ to
$[\frac{1}{3}a_{\varphi(i)}+\frac{2}{3}b_{\varphi(i)};\frac{2}{3}a_i+\frac{1}{3}b_i]$ ; 
\\
on $[b'_i;a'_{\psi(i)}]$, $h_i$ coincides with the affine map
from $[\frac{1}{3}a'_i+\frac{2}{3}b'_i;\frac{2}{3}a'_{\psi(i)}+\frac{1}{3}b'_{\psi(i)}]$ to
$[\frac{1}{3}a_i+\frac{2}{3}b_i;\frac{2}{3}a_i{\psi(i)}+\frac{1}{3}b_{\psi(i)}]$ ; 
\\ $h_i$ coincides with $h_{i-1}$ on 
$[0;\frac{1}{3}a'_{\varphi(i)}+\frac{2}{3}b'_{\varphi(i)}]\cup
[\frac{2}{3}a'_{\psi(i)}+\frac{1}{3}b'_{\psi(i)};1]$.
\end{itemize}

For all $i\in\N^*$ and $t\in[0;1[$, one defines then $h_{t,i}=th_i+(1-t)h_{i-1}$.
One considers then a non-decreasing sequence $(t_n)_{n\in\N}$ with $t_0=0$ and
$t_n\limite{n\to\infty}{}1$, and, given $t\in[0;1[$, if $t_i\leqslant t<t_{i+1}$, then
one defines $h_t=h_{\frac{t-t_i}{t_{i+1}-t_i},i+1}$.
One can check that the path of $\mc{C}^1$-diffeomorphisms constructed by this way is continuous.

Furthermore, one has : $h_tf_0h_t^{-1}\limite{t\to 1}{}f_1$ with respect to the $\mc{C}^1$-topology.

~~Indeed, let $\varepsilon$ be a strictly positive real number.
From the regularity of the diffeomorphisms $f_0$ and $f_1$, and since the fixed points of these two diffeomorphisms, except $0$ and $1$, are not hyperbolic, there exists a whole number $n\in\N$ such that, for all $i\geqslant n$, one has :
$\Vert\restriction{f_0}{C_i}-id\Vert_{1}<\dfrac{\varepsilon}{2}$ and
$\Vert\restriction{f_1}{C_i}-id\Vert_{1}<\dfrac{\varepsilon}{2}$.
In particular, if $t\geqslant 0$, since $h_t$ is an affine map on the intervals where $f_0$ does not coincide with the identity, one has also 
$\Vert \restriction{h_tf_0h_t^{-1}}{h_t(C_i)}-id\Vert_{1}<\dfrac{\varepsilon}{2}$.

On the other hand, if $t\geqslant t_{n-1}$, then $h_tf_0h_t^{-1}$ coincides with $f_1$ on $\bigcup\limits_{k<n}C'_k$.
\\
On the complement of this set, $h_tf_0h_t^{-1}$ and $f_1$ either coincide with the identity (outside $\bigcup\limits_{k\geqslant n}h_t(C_k)$ as regards $h_tf_0h_t^{-1}$
and outside $\bigcup\limits_{k\geqslant n}C'_k$ as regards $f_1$), or, as we just saw it, are $\dfrac{\varepsilon}{2}$-close to the identity with respect to the $\mc{C}^1$-topology.
One can deduce from all this that $\Vert h_tf_0h_t^{-1}-f_1\Vert_{1}<\varepsilon$ for all 
$t\geqslant t_{n-1}$, and the proof is complete.
\hfill$\Box$\medskip

~~Let us denote by $g_{\varepsilon}$ the increasing $\mc{C}^1$-diffeomorphism coinciding with $g$ on 
$\bigcup\limits_{i\in J_{\varepsilon}}C'_i$ and with $id$ elsewhere, and let us define the $\mc{C}^1$-diffeomorphism $\tilde{g}_{\varepsilon}$ by :
$\tilde{g}_{\varepsilon}$ has no fixed point and has same sign as $\sigma'(C'_i)$ on $C'_i$, and
coincides with $id$ outside $\bigcup\limits_{i\in J_\varepsilon}C'_i$.
\\

~~By Corollary \ref{c.généralisation}, one knows then that there exists an isotopy by conjugacy from $f_1$ to $\tilde{g}_\varepsilon$, and then Theorem \ref{t.2}
ensures the existence of an isotopy by conjugacy from $\tilde{g}_\varepsilon$
to $g_\varepsilon$.
\\

~~Thus, by Proposition \ref{p2.généralisation}, there exists an isotopy by conjugacy from $f$ to $g_\varepsilon$ ; as a consequence there exists a conjugate of $f$ $\varepsilon$-close to $g_\varepsilon$, thus also $2\varepsilon$-close to $g$, and that is true for all $\varepsilon>0$ : the second item of the theorem is consequently proved.
\\

~~The first item can be proved following the same reasoning ; the only difference is that one does not "cancel the waves" of $g$ which are smaller than $\varepsilon$ before applying the described method. Consequently, one uses the one-to-one map which maps $C'$ to $C$, which in this case does exist without hypothesis of finiteness on the signature of $g$, and, by using the intermediate diffeomorphisms as in the case 2.(except $g_\varepsilon$), one obtains the existence of isotopies from the ones to the others, and finally one puts it end to end by using Proposition 
\ref{p2.généralisation}, so as to obtain the existence of an isotopy by conjugacy from $f$ to $g$.

\hfill$\Box$\medskip
\\

\section{Annex}

\subsection{Statement and structure of the proof}

~~We give here the proof of Lemma \ref{l1.transfochemin}. This Lemma can be proved by several ways ; the method expounded here has been gracefully suggested by C. Bonatti.
\\
~~First of all, we will work toward reducing the problem thanks to several remarks.

For convenience reasons, we will first transpose the problem to the following equivalent proposition, whose proof is the subject of this Annex : 

\begin{prop}\label{p.annexe}
Let $f$ be a $\mc{C}^1$-diffeomorphism of $[0;+\infty[$ without fixed point except $0$, such that $f$ is a contraction (i.e. $f-id<0$ sur $]0;+\infty[$).
Then there exist two $\mc{C}^1$-contractions $f_+$ and $f_-$ of $[0;+\infty[$ such that :
\begin{itemize}
\item $Df_+(0)=Df(0)=Df_-(0)$ ;
\item for all $x>0$, one has : $f_-(x)<f(x)<f_+(x)$ ;
\item for all $x>0$, there exists $n\in\N^*$ such that $f_-^n(x)<f^{n+1}(x)<f^n(x)<f^{^n-1}(x)<f_+^n(x)$.
\end{itemize}
\end{prop}

\underline{Remark :} If $f_+$ and $f_-$ are two $\mc{C}^1$-contractions of $[0;+\infty[$ such that $f_-<f<f_+$ and satisfying the conclusions of the proposition above only on a neighbourhood of $0$, then $f_+$ and $f_-$ satisfy the conclusions of the proposition on the whole half-line. Thus, it will be sufficient to construct these two contractions on a neighbourhood of $0$, and then it will not make any difficulty to extend them to contractions remaining respectively greater and smaller than $f$.
\\

Let us now thus present what will be the structure of the proof :

The first step consists in proving that, if $f$ is embeddable in a $\mc{C}^1$-flow, then
$f$ satisfies the conclusions of Proposition \ref{p.annexe}.

\begin{lem}\label{l1.annexe}
Let $f$ be a $\mc{C}^1$-contraction of $[0;+\infty[$ such that $f$ is the time one-map of a $\mc{C}^1$-vector field $X$ on $[0;+\infty[$.
Then there exist two $\mc{C}^1$-contractions $f_+$ and $f_-$ of $[0;+\infty[$ satisfying the three conclusions listed in Proposition \ref{p.annexe}.
\end{lem}

We will hence reduce now the initial problem by the following result :

\begin{lem}\label{l2.annexe}
If $f$ is a $\mc{C}^1$-contraction of $[0;+\infty[$, then there exist two $\mc{C}^1$-contractions $g_+,g_-$ of $[0;+\infty[$ and two $\mc{C}^1$-vector fields $X_+,X_-$ of $[0;+\infty[$ such that :
\begin{itemize}
\item $g_-<f<g_+$ on $[0;+\infty[$ ;
\item $g_+,g_-$ are the respective time one-maps of $X_+,X_-$ ;
\item $Dg_-(0)=Df(0)=Dg_+(0)$.
\end{itemize}
\end{lem}

~~Indeed, one can easily check that, by applying Lemma \ref{l1.annexe} to the contractions $g_+$ and $g_-$, one finds the wanted contractions $f_+$ and $f_-$ for $f$.
\\

\subsection{Proof of Lemma \ref{l1.annexe}}

~~Let us denote by $\varphi(t,x)$ the flow whose $f$ is the time one-map.
By considering $1$ as origin, one can write each $x\in]0;1]$ in the form $x=\varphi(t(x),1)$.
The map $x\longmapsto t(x)$ is then a positive, strictly decreasing map, and tends to $+\infty$ when $x$ tends to $0$. Moreover, it is differentiable and one can easily check that it has derivative $\dfrac{1}{X(x)}$ for all $x\in]0;1]$.

One defines then $f_+$ and $f_-$ on $]0;f(1)[$ by the following :
$$f_+(x)=\varphi(1-\dfrac{1}{t(x)},x) \mbox{ and }
f_-(x)=\varphi(1+\dfrac{1}{t(x)},x)$$

One can check that $f_+,f_-$ satisfy
$f_-<f<f_+<id$ on $]0;f(1)[$.
Let us show now that these two formulas define $\mc{C}^1$-contractions of $]0;f(1)[$ with same derivative as $f$ at $0$.

Indeed, one has the following calculation for $f_-$ :
$$\dfrac{d }{dx}(\varphi(1+\dfrac{1}{t(x)},x))=
X(\varphi(1+\dfrac{1}{t(x)},x))\cdot\dfrac{1}{X(x)}-
X(\varphi(1+\dfrac{1}{t(x)}))\cdot\dfrac{1}{X(x)}\cdot\dfrac{1}{t^2(x)}\mbox{,}$$ the first element coming from the calculation of the derivative of $\varphi$ with respect to its second variable and the second one from the derivative of $\varphi$ with respect to the time.

Notice that $\dfrac{d }{dx}(\varphi(1+\dfrac{1}{t(x)},x))=X(\varphi(1+\dfrac{1}{t(x)},x))\cdot\dfrac{1}{X(x)}\cdot(1-\dfrac1{t^2(x)})$ is strictly positive, as $t(x)>1$ on $]0,f(1)[$.  One deduces that $f_-$ is a diffeomorphism from $]0,f(1)[$ on its image. 

It remains to show that $f_-$ is $C^1$ at $0$ and that $Df_-(0)=Df(0)$. For that, we just have to notice that the first element has same limit when $x$ tends to $0$ as $\dfrac{X(f(x))}{x}$, and this expression tends to $Df(0)$ when $x$ tends to $0$.
The second element, by the same reasoning, tends to $0$.
The case of $f_+$ can be concluded by analogous calculations.

~~What remains to show now is that the third conclusion of Proposition \ref{p.annexe} is satisfied. One can write, for all $x\in]0;1]$ and for all whole number $n$, $f^n(x)=\varphi(n,x)$.
One denotes now $f_-^n(x)=x_n$ ; one has : $x_n=\varphi(t_n,x)$, where 
$t_n=n+\dfrac{1}{t(x_{n-1})}+\dfrac{1}{t(x_{n-2})}+\ldots+\dfrac{1}{t(x_1)}+\dfrac{1}{t(x)}$. If we show now that this sum $\dfrac{1}{t(x_{n-1})}+\dfrac{1}{t(x_{n-2})}+\ldots+\dfrac{1}{t(x_1)}+\dfrac{1}{t(x)}$ is greater than $1$, then the proof will be complete.
We do it by contradiction : if this sum were bounded by some constant $C$, then each
$t(x_n)$ should be smaller than $n+C$, hence the sum $\dfrac{1}{t(x_{n-1})}+\dfrac{1}{t(x_{n-2})}+\ldots+\dfrac{1}{t(x_1)}+\dfrac{1}{t(x)}$ should be greater than $\sum\limits_{k=1}^n\dfrac{1}{k+C}$, which can be arbitrarily great. Thus we have a contradiction.

One can finally extend the constructed diffeomorphisms $f-+$ and $f_-$ restricted to a small neighbourhood of $0$ to $\mc{C}^1$-contractions of $[0;+\infty[$ in such a way that the inequality $f_-<f<f_+$ is preserved. 
\hfill$\Box$\medskip

The matter is thus now to prove Lemma \ref{l2.annexe}, and then the proof will be complete.

\subsection{Proof of Lemma \ref{l2.annexe} :}

~~In the case where $Df(0)=1$, one has the following proof :
\\
One defines the two following $\mc{C}^1$-vector fields on $[0;1]$ :
$$X_-(x)=\dfrac{1}{2}(f(x)-x)\dfrac{d}{dx},X_+(x)=2(f(x)-x)\dfrac{d}{dx}\mbox{ .}$$
Then, one can check that these two vector fields satisfy the wanted properties, by using the fact that $f$ is $\mc{C}^1$-close to the identity on a neighbourhood of $0$.
\\

~~The case where $Df(0)\neq 1$ is quite more difficult, and in order to solve it two steps will be needed.

The first one consists in noticing that, if the derivative of $f$ is a monotonic map, then $f$ satisfies the conclusions of Lemma \ref{l2.annexe} :

\begin{lem}\label{l3.annexe}
Let $f$ be a $\mc{C}^1$-contraction of $[0;+\infty[$ such that its derivative is a monotonic continuous map which is not equal to $1$ at $0$. Then there exist two $\mc{C}^1$-contractions $g_+,g_-$ of $[0;+\infty[$ and two $\mc{C}^1$-vector fields $X_+,X_-$ of $[0;+\infty[$ satisfying the two conclusions listed in Lemma \ref{l2.annexe}.
\end{lem}

Finally, the following lemma will enable us to complete the proof :

\begin{lem}\label{l4.annexe}
Let $f$ be a $\mc{C}^1$-contraction of $[0;+\infty[$ such that $Df(0)\neq 1$.
Then there exist two $\mc{C}^1$-contractions $h_+,h_-$ of $[0;+\infty[$ such that :
\begin{itemize}
\item $h_-<f<h_+$ on some interval $]0;\varepsilon]$, where $\varepsilon>0$ ;
\item $Dh_-(0)=Df(0)=Dh_+(0)$ ;
\item $h_-$ and $h_+$ have monotonic derivative.
\end{itemize}
\end{lem}

Indeed, given two such contractions $h_+,h_-$ for a contraction $f$, one can apply Lemma \ref{l3.annexe} to each of these two maps in order to obtain two $\mc{C}^1$-contractions which will be the time one-maps of $\mc{C}^1$-vector fields on $[0;+\infty[$ and be respectively  greater and smaller than $f$.

\underline{Proof of Lemma \ref{l3.annexe} :}
~~We begin by showing this lemma inthe case where $Df$ is increasing.
Since $f$ is greater than the homothety $x\longmapsto Df(0)\cdot x$, one can define 
$g_-(x)=Df(0)\cdot x$, which is the time one-map of the linear vector field 
$X_-(x)=\log(Df(0))\cdot x \dfrac{d}{dx}$ and satisfies the wanted conditions.

As regards $g_+$, we will consider the vector field of $]0;f(1)[$ given by :
$X_+(x)=\log(\dfrac{x}{f^{-1}(x)})\cdot x \dfrac{d}{dx}$, and we define $g_+$ as being the time one-map of this vector field.
As first remark, let us notice that, since $Df$ is increasing, one has also that $\dfrac{f(x)}{x}$ is increasing ; then 
$\vert\log(\dfrac{y}{f^{-1}(y)})\vert\leqslant\vert\log(\dfrac{f(x)}{x})\vert$ for all
$y\in[f(x);x]$.
Given $x_0\in ]0;f(1)[$, for all $y\in [f(x_0);x_0]$, the norm of the vector field $X^+(y)$ is smaller than the one of $\log(\dfrac{f(x_0)}{x_0})\cdot x \dfrac{d}{dx}$, whose time one-map maps $x_0$ to $f(x_0)$.
It follows that it takes to the point $x_0$ a time longer than $1$ to reach $f(x_0)$ along the orbit of $X_+$. Hence $g_+(x_0)\in[f(x_0);x_0]$, and thus $g_+(x_0)\geqslant f(x_0)$ : the map $g_+$ is greater than $f$.

One has now to check that $g_+$ has $\mc{C}^1$-regularity also at $0$, and that its derivative at this point is the same as $f$ : it follows from the expression of the derivative of $X_+$, and then it suffices to make tend $x$ to $0$ in this expression.

~~In the case where $Df$ is decreasing, we follow the same reasoning with 
$g_+(x)=Df(0)\cdot x$ and $g_-$ being the time one-map of the vector field
$X_-(x)=\log(\dfrac{x}{f^{-1}(x)})\cdot x\dfrac{d}{dx}$.

\hfill$\Box$\medskip

\underline{Proof of Lemma \ref{l4.annexe} :}
One defines $h_+$ and $h_-$ on a neighbourhood of $0$ by :

$h_+(x)=\int\limits_{0}^{x}\inf\limits_{y\in[0;z]}Df(y)\mathrm{dz}$ and
$h_-(x)=\int\limits_{0}^{x}\sup\limits_{y\in[0;z]}Df(y)\mathrm{dz}$.
These two expressions define $\mc{C}^1$-maps whose derivatives are strictly smaller than $1$, thus they define on a neighbourhood of $0$ $\mc{C}^1$-contractions. Their derivative at $0$ is $Df(0)$. All what is needed now is to  extend $h_+$ and $h_-$ as
$\mc{C}^1$-diffeomorphisms on the half-line.

\section{Appendix}

We present here the argument given by A. Navas, which proves Theorem~\ref{t.navas}. The proof is based on ergodic theory and reduces the problem to the approximate solving of a  cohomological equation. The starting point is the following remark : 

Let $f$ be a $\mc{C}^1$-diffeomorphism of the interval without hyperbolic fixed point.
We want to find a $\mc{C}^1$-diffeomorphism $\varphi$ of $[0;1]$ which conjugates $f$ to a $\mc{C}^1$-diffeomorphism which would be close to the identity ; it amounts to finding $\varphi$ such that 
$\log D(\varphi f\varphi^{-1})$ is close to $0$, in other words such that
$\log(D\varphi)-\log(D\varphi)\circ f$ is close to $\log Df$. Then we want to find approximate and continuous solutions to the cohomological equation 
\[ \rho-\rho\circ f=\log Df \mbox{.}\]

By if necessary adding to $\rho$ a constant, the map $\varphi(x)=\int_{0}^{x} \exp(\rho(u))du$ becomes a $\mc{C}^1$-diffeomorphism of $[0;1]$ satisfying the wanted conditions. More precisely, that argument shows : 

\begin{lem} Let $f$ be a diffeomorphism of the interval or of the circle $S^1=\R/\Z$.  Let $(\rho_t)_{t\in[0,1[}$ be a continuous path of continuous maps such that $\rho_t-\rho_t\circ f$ converges uniformly to $\log Df$ when $t$ tends to $1$ and such that $\int_0^1 \exp{\rho_t(u)}du=1$. One defines 
 $$h_t(x)=\int_{0}^{x} \exp(\rho_t(u))du.$$

Then $(h_t)_{t\in[0,1[}$ is a continuous path of $\mc{C}^1$-diffeo\-mor\-phisms of the interval (or of the circle) and form an isotopy by conjugacy from $f$ to the identity (in the case of the interval) or to the rotation of same rotation number as $f$ (in the case of the circle). 
\end{lem}

Then we have to prove the existence of the path $\rho_t$. The following proposition ensures the existence of approximate solutions to the cohomological equation.

\begin{prop}\label{p.appendice}
Let $f$ be a $\mc{C}^1$-diffeomorphism of $[0;1]$ without hyperbolic fixed point (resp. a diffeomorphism of the circle with irrational number of rotation), and $\varepsilon>0$ be a real number.
Then there exists a continuous map $\rho$ of $[0;1]$ (resp. of $S^1$) such that $\Vert \rho-\rho\circ f-\log(Df)\Vert_{\infty}<\varepsilon$.
\end{prop}

\underline{Proof :} 

Given $k\geqslant 1$ a whole number and $\varphi$ a map from $[0;1]$ to $\R$, one considers the $k$-th Birkhoff sum of $S_k(\varphi)=\sum_{i=0}^{k-1}\varphi(f^i)$, 
and, given $n$ a whole number, one considers the map $\rho_n$ defined as follows :
\\
$\rho_n=\dfrac{S_1(\log(Df))+S_2(\log(Df))+\ldots +S_n(\log(Df))}{n}$.
One can then calculate 
$\rho_n-\rho_{n+1}\circ f=\log(Df)-\dfrac{S_n(\log(Df))\circ f}{n}$. One concludes by the following claim : 

\begin{affi}
 $\dfrac{S_n(\log(Df))}{n}$ converges uniformly to $0$ when $n$ tends to $\infty$.
\end{affi}

\underline{Proof of the Claim :} This is a consequence of the fact that the hypotheses \emph{"no hyperbolic fixed point"} and \emph{"irrational number of rotation"} imply that, for all probability measure $\mu$ which is invariant under $f$ and ergodic, the Lyapunov exponent of $f$ with respect to $\mu$, given by $\lambda(\mu)=\int \log (Df)d\mu$, is equal to $0$. See \cite{PS} or \cite{H} to have more details about this result. 

\hfill$\Box$\medskip

~~Now we are able to construct a sequence
$(h_n=\int_{0}^{x}\exp(\rho(u))du)_{n\in\N}$ of conjugacies such that
$h_n f h_n^{-1}$ converges to the identity with respect to the $\mc{C}^1$-topology.
We have still to obtain, not only a sequence of conjugates of $f$, but an isotopy by conjugacy from $f$ to the identity.

Given a whole number $n$ and a real number $\lambda$ between $0$ and $1$, if $t$ is equal to $\dfrac{\lambda}{n}+\dfrac{1-\lambda}{n+1}$, we will denote by $\rho_t$ the map (of the interval or of the circle) defined by  :
$$\rho_t=\lambda\rho_n+(1-\lambda)\rho_{n+1}.$$

One can now conclude by the following lemma, whose proof is immediate :
\begin{lem} The path $(\rho_t)_{t\in[0,1[}$ is a continuous path of continuous maps such that, for all $t$, one has $\int_0^1 \exp{\rho_t(u)}du=1$ and such that $\rho_t-\rho_t\circ f$ converges uniformly to $\log Df$ when $t$ tends to $1$.
\end{lem}


\vskip 5mm
\begin{tabular}{l}
\'Eglantine Farinelli
\\
\footnotesize{eglantine.farinelli@u-bourgogne.fr}
\\
Institut de Math\'ematiques de Bourgogne
\\
CNRS - URM 5584
\\
Universit\'e de Bourgogne
\\
Dijon 21000, France

\end{tabular}

\end{document}